\documentclass[9pt,leqno]{smfart}

\author{Fabrice Orgogozo}
\address{Princeton University, Mathematics Department\\ 
Fine Hall, Washington Road\\Princeton, NJ~\osn{8544}---\osn{1000}
\ U.S.A.}
\email{orgogozo@math.princeton.edu}

\newcommand{\timesf}{{\sr{\la}{\times}}}

\usepackage{amsmath,a4wide}
\usepackage{amssymb}
\usepackage[latin1]{inputenc}
\usepackage[T1]{fontenc}
\usepackage[all]{xy}
\usepackage[francais]{babel}
\usepackage{aeguill}
\usepackage{smfthm}
\def«{\og}
\def»{\fg}

\SwapTheoremNumbers

\makeatletter
\renewcommand\theequation{\thesubsection.\textbf{\alph{equation}}}
\@addtoreset{equation}{subsection}
\makeatother

\renewcommand\thesubsubsection{\textbf{\thesubsection.\arabic{subsubsection}}} 
\renewcommand\thesubsection{\textbf{\thesection.\arabic{subsection}}}

\newtheorem{prp}[subsection]{Proposition} 
\newtheorem{prp2}[subsubsection]{Proposition}
\newtheorem{thm}[subsection]{Théorème}

\newtheorem{lmm}[subsection]{Lemme}
\newtheorem{lmm2}[subsubsection]{Lemme}

\newtheorem{crl}[subsection]{Corollaire}

\theoremstyle{definition}

\newtheorem{rmr}[subsection]{Remarque}
\newtheorem{rmr2}[subsubsection]{Remarque}

\newtheorem{rmrs}[subsection]{Remarques}

\def\sga#1#2#3{[{\bf $\mathbf{SGA_{#1}}$}~{\sc #2}~#3]}
\def\ega#1#2{[{\bf ÉGA}~{\sc #1}~#2]}
\newcommand{\beqt}{\begin{equation}}
\newcommand{\eeqt}{\end{equation}}

\newcommand{\Hom}{\mathsf{Hom}}
\newcommand{\Aut}{\mathsf{Aut}}

\newcommand{\ob}{\mathsf{Ob}\ }

\newcommand{\Ens}{\mathsf{Ens}}

\newcommand{\Sch}{\mathsf{Sch}}

\newcommand{\SP}{\mathrm{Spec}}
\newcommand{\aff}{\mathbf{A}}
\newcommand{\etb}{\bar{\eta}}

\newcommand{\red}{\text{réd}}

\newcommand{\cp}{\Psi}
\newcommand{\ce}{\Phi}
\newcommand{\RG}{\mathrm{R\Gamma}}

\newcommand{\HH}{\mathrm{H}}

\newcommand{\gp}{\pi_1}     

\newcommand{\dbc}{\mathsf{D}_c^b}

\newcommand{\ZZ}{\mathbf{Z}}

\newcommand{\NN}{\mathbf{N}}

\newcommand{\QQ}{\mathbf{Q}}

\newcommand{\CC}{\mathbf{C}}
\newcommand{\PP}{\mathbf{P}}
\newcommand{\FF}{\mathbf{F}}

\newcommand{\QL}{{\QQ_{\ell}}}

\newcommand{\ra}{\rightarrow}
\newcommand{\la}{\leftarrow}
\newcommand{\Rra}{\Rightarrow}
\newcommand{\rra}{\rightrightarrows}
\newcommand{\lra}{\longrightarrow} 

\newcommand{\hra}{\hookrightarrow }
\newcommand{\sr}{\stackrel}

\newcommand{\rrraxy}{\ar@<1ex>[r] \ar@<-1ex>[r] \ar[r] }
\newcommand{\iso}{\stackrel{\sim}{\ra}}
\newcommand{\giso}{\stackrel{\sim}{\la}}
\newcommand{\isononcan}{\simeq}

\def\fl#1{\overleftarrow{#1}}

\newcommand{\R}{\mathrm{R}}

\def\]{\textup{\mbox{]\hspace{-.15em}]}}}
\def\[{\textup{\mbox{[\hspace{-.15em}[}}}

\def\mc{\mathcal}
\def\got{\mathbf}

\def\-l{\ \vspace{-4mm}}
\def\trdist#1#2#3{ \xymatrix{#1 \ar[rr] & & #2 \ar[dl] \\ & #3 \ar[ul]^{+1} & }}
\def\L#1{\mathsf{L}_{#1}}
\newcommand{\cad}{c'est-à-dire\xspace}

\def\sur{\overline}

\def\osn{\oldstylenums}

\title{Modifications et cycles évanescents sur une base de dimension
supérieure à un.}

\thanks{Ce travail a été effectué à l'université de Princeton (É.-U.A) ; l'auteur
a bénéficié de l'excellent accueil de son département de mathématique et en particulier
de son directeur de l'époque, Nick Katz.}

\begin{document}

\maketitle \date{18 juillet 2005}
\section*{Introduction}

On se propose de démontrer l'analogue en cohomologie étale du théorème
principal de \cite{cycles_evanescents@Sabbah}, conjecturé par Pierre
Deligne, selon lequel tout morphisme de schémas acquiert un
comportement de «morphisme sans éclatement»,  quitte à effectuer une
modification de la base.  Nous renvoyons à \ref{énoncé} pour un énoncé
précis en termes de cycles évanescents (ou plutôt «proches») ---~dont
la définition est rappelée dans la section \S 1~--- sur une base de dimension
quelconque. Pour un exemple d'éclatement «avec éclatement» (sic !), voir \S \ref{unexemple}.

Le théorème principal entraîne en particulier que sur une telle base,
si les cycles évanescents ne commutent plus nécessairement aux
changements de base (même si celle-ci est régulière),  cela est malgré
tout vrai  \emph{localement pour la topologie de la descente effective
universelle}.

Ces résultats généralisent également ceux de Roland Huber
\cite{Adic_spaces@Huber}, \S 4.2, sur la cohomologie des cycles
proches quand la base est le spectre  d'un anneau valuatif
(cf. \ref{remarque-adique}). Ces derniers sont une étape
importante dans la démonstration des théorèmes de commutation au
changement de base et constructibilité de \emph{op. cit.}

Les «cycles évanescents» dont il est question ont été introduits par
P.~Deligne à l'occasion d'un séminaire à l'IHÉS sur les fonctions $L$
il y a une vingtaine d'années et généralisent ceux définis par
Alexandre Grothendieck dans \sga{7}{i}{}\nocite{sga7}.  C'est un moyen
commode de regrouper la cohomologie des différents «tubes de Milnor»
\mbox{$X(x)\times_{S(s)} S(t)$} où $t$ est une générisation de
$s=f(x)$ (les points sont géométriques), lorsque que l'on  se donne un
morphisme de schémas $X\sr{f}{\ra} S$ et un faisceau étale $\mc{F}$
sur $X$.  Ces tubes apparaissent
naturellement quand on calcule la fibre en $t$ d'un faisceau $\R^i
{f_x}_*\mc{F}$,  où $f_x$ est le morphisme local $X(x)\ra S(s)$ induit
par $f$ entre les hensélisés stricts en $x$ et $s$ (les «petites
boules» centrées en ces points). Dans ce langage, le théorème
principal affirme  qu'après modification de la base, les morphismes
$f_x$ se comportent  cohomologiquement comme des morphismes propres.\\
L'ingrédient essentiel de notre démonstration est un théorème sur les
fibrations plurinodales dû à A.~Johan de~Jong (\cite{Families@de_Jong}).

Enfin, en guise d'application de ces topos exotiques, mais dans le cas 
beaucoup plus simple des singularités isolées, nous donnons une démonstration
de la conjugaison des cycles évanescents associés à un pinceau de Lefschetz,
plus proche de l'intuition topologique et également valable en caractéristique
deux. L'ingrédient clé est ici dû à Ofer Gabber.

Je remercie très chaleureusement \textsc{Ofer Gabber} et \textsc{Luc
Illusie} de m'avoir suggéré ce problème ainsi que d'en avoir discuté
avec moi.  C'est également avec plaisir que je remercie \textsc{János
Kollár} d'avoir gentiment répondu  à quelques questions naïves et le
rapporteur pour sa relecture extrêmement attentive et ses commentaires utiles et détaillés.
Enfin, je suis spécialement reconnaissant envers \textsc{Ofer Gabber}
de m'avoir expliqué la proposition \ref{lemme-Ofer} à un moment
critique, ainsi que de m'avoir fait part de ses précieuses remarques à l'origine 
de nombreuses améliorations.

\part{Commutation aux changements de base}

\section{Rappels sur les cycles évanescents sur une base de dimension
supérieure à $1$ et énoncé du théorème}\label{rappels}

On reprend brièvement les définitions de \cite{Vanishing@Laumon},
auquel on renvoie le lecteur pour plus de détails.  Pour tout
morphisme de topos $\got{f}:\got{X}\ra \got{S}$, on peut définir un
troisième topos $\got{X}\timesf_{\got{S}} \got{S}$
(\textit{loc. cit.}, \S~3 ou  \cite{Proper@Moerdijk},
Chap. \textsc{v}, \S~1) solution d'un problème  universel
$2$-catégorique, dont les points sont les paires $(\got{x},\got{t})$
où $\got{x}$ est un point de $\got{X}$ et  $\got{t}$ une générisation
de $\got{f}(\got{x})$.  Ce topos est naturellement équipé d'un
morphisme $\Psi_{\got{f}}:\got{X} \ra  \got{X}\timesf_{\got{S}}
\got{S}$, correspondant sur les points à $\got{x}\mapsto
(\got{x},\got{f}(\got{x}) \sr{\mathrm{id}}{\ra} \got{f}(\got{x}))$
ainsi que de deux projections $p:\got{X}\timesf_{\got{S}} \got{S} \ra
\got{X}$ et $q:\got{X}\timesf_{\got{S}} \got{S} \ra  \got{S}$.

Si l'on suppose maintenant que $\got{X}$ (resp. $\got{S}$)  est le
topos étale associé à un schéma $X$ (resp. $S$), et $\got{f}$ induit
par un morphisme de schémas $f:X\ra S$, on montre que pour tout
faisceau $\mc{F}$ sur $\got{X}$, la fibre de
$\R\Psi_{\got{f}\,*}\mc{F}$  en un point $(\got{x},\got{t})$ est
isomorphe à $\RG(X(\got{x})\times_{S(\got{f(x)})} S(\got{t}),\mc{F})$.

Dorénavant, on confondra $X$ avec le topos étale associé, afin
d'alléger les notations.

Un anneau commutatif $\Lambda$ étant donné, le foncteur «cycles
proches» de $f$ est le foncteur dérivé  $${\R\Psi}_{f\,*}^{\Lambda} :
\mathsf{D}^{+}(X,\Lambda)\ra  \mathsf{D}^{+}(X\timesf_S S,\Lambda).$$
On le notera généralement $\cp_{f}$, l'anneau $\Lambda$ étant
sous-entendu.

Avant d'énoncer le théorème principal rappelons qu'une
\emph{modification} (resp.  une \emph{altération}) $S'\ra S$ est un
morphisme propre surjectif  induisant un isomorphisme (resp. un
morphisme fini) au-dessus d'un ouvert partout dense de $S$, tel que
tout point  maximal de $S^{\prime}$ s'envoie sur un point maximal de
$S$.  Enfin on rappelle que l'on abrège «quasi-compact et
quasi-séparé» en \textit{cohérent}.

\begin{thm}\label{énoncé} Soient $S$ un schéma cohérent irréductible,
$f:X\ra S$ un morphisme de présentation finie, $n$ un entier
inversible sur $S$ et $\mc{F}\in \ob \dbc(X,\ZZ/n\ZZ)$.  Alors, le
complexe ${\cp_{f}}(\mc{F})$ n'a pas de cohomologie en grand degré et
il existe une \emph{modification} $S'\ra S$ telle que l'on ait la
propriété suivante : notant $f_{S'}$ (resp. $\mc{F}_{S'}$) le
morphisme $f\times_S S':X'\ra S'$ (resp. le complexe $\mc{F}_{|X'}$),
la \emph{formation de $\cp_{f_{S'}}\mc{F}_{S'}$ commute à tous les
changements de base $g:T\ra S'$}.\\ Plus précisément, le morphisme  de
changement de base ([\textbf{SGA} 4] \textsc{xii}.4 et \textsc{xvii}.2.1)

\begin{equation}\label{cb} c_{f,\mc{F},T/S'} :
\fl{g}^*\big(\cp_{f_{S'}}\mc{F}_{S'}\big) \ra \cp_{f_{T}}\mc{F}_{T}
\end{equation}

associé au diagramme essentiellement commutatif de topos :

$$
\xymatrix{ X_{S'} \ar[d]^{\cp_{f_{S'}}} & X_{T} \ar[l]^g
\ar[d]^{\cp_{f_{T}}}  \\ X_{S'}\timesf_{S'} S' & X_T \timesf_{T} T
\ar[l]^{\fl{g}} } $$

est un \emph{isomorphisme}.
\end{thm}

Nous verrons dans la seconde partie que l'on peut améliorer cet énoncé
(comparaison entre la cohomologie des tubes de Milnor et  celle des
fibres de Milnor~(\ref{tubes-fibres}) ; constructibilité
(\ref{constructibilite})).

Insistons sur le fait que le morphisme $S'/S$ dépend de $\mc{F}$.

La démonstration se fait par une triple récurrence sur un indice de troncation, la dimension
de $S$ (que l'on peut supposer finie comme expliqué plus bas) et  la
«dimension relative générique» de $f$. On procède par réduction au cas
des courbes semi-stables et des faisceaux constants.  L'entier $n$
étant choisi une fois pour toute, on écrira souvent $\Lambda$ pour $\ZZ/n\ZZ$.

Le morphisme $f:X\ra S$ étant sous-entendu,  pour chaque point $(x,t)$
de $X\timesf_S S$, c'est-à-dire la donnée un point géométrique $x$ de
$X$ et une générisation géométrique $t$ du point géométrique image de
$x$ dans $S$, on désignera par $X_{x,t}$ le schéma produit fibré
\mbox{$X(x)\times_{S(f(x))} S(t)$} ; c'est le \textit{tube de Milnor} (de
$f$) en $(x,t)$.

La démonstration de \ref{énoncé} occupe les trois sections suivantes.
Le schéma de la démonstration est expliqué dans le
paragraphe~\ref{devissage-ultime}, qui conclut les dévissages
liminaires et détaille la récurrence.

Notation : pour $m$ un morphisme de topos, on écrira principalement $m_*$ pour l'image
directe dérivée $\R m_*$.

\section{Premières réductions}\label{réductions}

\subsection{Passage à la limite}\label{limite}

Soit $(S_i)_{i\in I}$ un recouvrement ouvert affine de $S$ ; par
quasi-compacité de $S$, on peut supposer $I$ fini. Procédant comme en
\cite{Platification@Raynaud}, première partie, \S~5.3.3, on voit qu'il
suffit de démontrer le théorème pour $S$ affine.  Le morphisme $f$
étant de présentation finie, il existe d'après \ega{iv}{8.9.1},  un
$\SP(\ZZ)$-schéma affine de type fini $S_0$, un morphisme de type fini
$f_0 : X_0 \ra S_0$, un morphisme dominant $S\ra S_0$ et un
$S$-isomorphisme $X_0\times_{S_0} S\iso X$.\\ De même, on peut
supposer que le complexe borné $\Lambda$-\emph{constructible} $\mc{F}$  est
isomorphe à l'image inverse d'un tel complexe sur $X_0$.

Pour cette raison, les dévissages qui vont suivre concernent une base
noethérienne.

\subsection{Un résultat d'annulation}

Le fait que le complexe $\cp_{f}\mc{F}$, \textit{a priori} dans
$\mathsf{D}^+(X\timesf_S S,\Lambda)$, soit dans
$\mathsf{D}^b(X\timesf_S S,\Lambda)$ résulte de la proposition
suivante :

\begin{prp}\label{annulation} Soient $N$ un entier, $S$ un schéma et
$f:X\ra S$ un morphisme de type fini,  dont les fibres sont de
dimension inférieure ou égale à $N$.  Si $\mc{F}$ est un faisceau abélien
de torsion sur $X$,  pour tout point $(x,t)$ de $X\timesf_S S$,
$$\HH^i{\cp_f\mc{F}}_{(x,t)} \isononcan \HH^i(X(x)\times_{S(f(x))}
S(t),\mc{F})=0$$  dès que $i>2N$.
\end{prp}

Faute de référence, voici une  démonstration, qui m'a été communiquée
par O.~Gabber, évitant ainsi d'avoir recours à une variante tronquée
des énoncés principaux.

\begin{proof} La question étant locale au voisinage de chaque point de
$X$, on peut supposer $X$ affine et même quasi-fini sur le $S$-espace
affine $\aff^N_S$, avec $S$ également affine. Appliquant le Main
Theorem de Zariski et  prolongeant le faisceau $\mc{F}$ par zéro,  on
peut supposer $X$ fini sur $\aff^n_S$ puis finalement $X/S$
\textit{projectif} à fibres de dimension inférieure à $N$.  Soit
$(x,t)$ un couple de points géométriques comme dans l'énoncé ; notons
$s=f(x)$ et supposons pour simplifier les notations que $S(s)\iso S$.
Notons $\tau$ le morphisme essentiellement étale $X(x)\times_{S} S(t)
\ra X \times_{S} S(t)$. On a  $\RG(X(x)\times_{S}
S(t),\mc{F})\isononcan\RG(X\times_{S} S(t), \R\tau_*\mc{F})$. De plus,
d'après le théorème de changement de base propre,  $\RG(X\times_{S}
S(t),\R\tau_*\mc{F})\iso\RG(X_t,(\R\tau_*\mc{F})_{X_t})$.  Comme $X_t$
est de $n$-dimension cohomologique étale $\leq 2N$, il suffit de
montrer que $\HH^i((\R\tau_*\mc{F})_{X_t})=0$ pour $i>0$.  Cela
signifie que pour tout point géométrique $u$ de $X_t$,
$(\R^i\tau_*\mc{F})_u=0$ pour $i>0$ (on identifie $u$ et son image par
le morphisme $X_t\ra X\times_{S} S(t)$). Cette fibre est isomorphe à
$\HH^i(X(x)\times_X X(u),\mc{F})$ ; d'après \cite{Joins@Artin},~3.4,
les composantes connexes de $X(x)\times_X X(u)$  sont strictement
locales. Le résultat en découle. L'hypothèse de projectivité est là
pour nous assurer que $u$ et $x$ sont contenus dans un même ouvert
affine.
\end{proof}

\subsection{Réduction au cas  d'un faisceau constant}\label{faisceau constant}

Notons que l'on peut supposer dans la démonstration de~\ref{énoncé}
que $\mc{F}$ est un faisceau concentré en degré $0$.  Comme $\mc{F}$
est constructible, il existe un nombre fini de morphismes finis
$(p_i:X_i\ra X)_{i\in I}$,  tels que $\mc{F}$ s'injecte dans le
faisceau  $\displaystyle \prod_{i\in I} p_{i*}C_i$, où chaque $C_i$ est
un  faisceau \textit{constant} constructible de $\Lambda$-modules sur $X_i$.
Par commodité nous dirons dans ce paragraphe que ces faisceaux sont
\emph{très constructibles} sur $X$.  Supposons \ref{énoncé} démontré
pour les faisceaux très constructibles.  Remarquons que l'image
inverse d'un faisceau très constructible est très constructible
(changement de base pour les morphismes finis).  Une récurrence immédiate sur
les entiers $k$ tels qu'après modification les cônes des morphismes
(\ref{cb}) n'aient pas de cohomologie en degré inférieur ou égal à $k$, montre qu'il
suffit, pour un faisceau $\mc{F}$ donné, de démontrer \ref{énoncé}
pour  un nombre fini de faisceaux très constructibles (cogénérateurs
de la catégorie des faisceaux constructibles).  Cela découle de la
proposition précédente car il suffit donc de
démontrer le résultat pour un nombre fini d'entiers $k$.  Ainsi, on
s'est ramené au cas où $\mc{F}=p_*C$, pour un morphisme fini $p:X'\ra
X$ et un faisceau $C$ constant $\Lambda$-constructible sur $X'$.  Si
la commutation aux changements de base est démontrée pour $(X'/S,C)$,
elle l'est aussi pour $(X/S,p_*C)$. En effet, pour tout point $(x,t)$
de $X\timesf_S S$,  $$ \RG(X_{x,t},p_*C)=\prod_{x_i\mapsto x}
\RG(X'_{x_i,t},C), $$ et ce universellement sur $S$, puisque $p_*C$
commute aux changements de base.\\ Il est utile de remarquer pour la
suite de la démonstration que quand  $S$ est irréductible, passer de
$X/S$ à $X'/S$ n'augmente pas la dimension de la fibre générique.

\subsection{Des altérations aux modifications}\label{altmod} Il suffit
de montrer qu'il existe une \textit{altération} $S'\ra S$ telle que
les conclusions de \ref{énoncé} soient satisfaites.  Cela résulte
immédiatement des deux lemmes suivants.

\begin{lmm}\label{platification} Soit $S'\ra S$ une altération entre
schémas noethériens.  Il existe deux modifications $\tilde{S}\ra S'$,
$S''\ra S$ et un morphisme fini surjectif $\tilde{S}\ra S''$ tels que
le  carré  $$ \xymatrix{ S' \ar[d]^{\mathrm{alt.}} & \tilde{S}
\ar@{->>}[d]^{\mathrm{fini}} \ar[l]^{\mathrm{modif.}} \\ S & S''
\ar[l]^{\mathrm{modif.}}  } $$ soit commutatif.
\end{lmm}

\begin{proof} C'est un corollaire immédiat de
\cite{Platification@Raynaud}, \textsc{i}~5.2.2, compte tenu du fait
qu'un morphisme propre, plat, génériquement fini est fini et que le
transformé propre de $S'/S$ via $S''\ra S$ domine $S'$.
\end{proof}

Je remercie O.~Gabber de m'avoir appris ce fait.

\begin{lmm}\label{descente-fini} Sous les hypothèses
\emph{\ref{énoncé}}, s'il existe  $S'/S$ \textrm{fini} surjectif tel
que la formation de $\cp_{f_{S'}}\mc{F}'$ commute aux changements de
base $T'\ra S'$, alors la formation du complexe $\cp_f \mc{F}$ commute
à tous les changements de base $T\ra S$.
\end{lmm}

\begin{proof} Soient $T\ra S$ un morphisme et $(x_T,b_T)$ un point de
$X_T \timesf_{T} T$ s'envoyant sur le point $(x,b)$ de $X\timesf_S S$.
Il s'agit de montrer que le morphisme canonique
$\RG(X_{x,b},\mc{F})\ra \RG({X_T}_{x_T,b_T},\mc{F})$ est un
isomorphisme, sachant que cette propriété est vraie après changement
de base par $S'/S$, fini surjectif.  On peut supposer $T$ et $S$
strictement locaux, de points fermés les images de $x$ et $x_T$.  Pour
chaque entier $q\geq 0$, notons  $S^{(q)}=(S'/S)^{q+1}$ le produit
fibré itéré (resp. $T^{(q)}=(T\times_S S'/T)^{q+1}$,
$X^{(q)}=X\times_S S^{(q)}$, $X_T^{(q)}=X_T\times_T T^{(q)}$).  Ils
définissent naturellement des hyperrecouvrements
propres\footnote{Rappelons brièvement (cf.
\cite{Hodge3@Deligne}, \S~5)  qu'un \textit{hyperrecouvrement propre}
d'un schéma $S$ est un objet simplicial $S_{\bullet}$ de la catégorie
des $S$-schémas (\textit{i.e.} un foncteur $\Delta^{\textrm{op}}\ra
\Sch_S$) tel que pour chaque entier $n\geq 0$,  le morphisme canonique
$S_{n+1}\ra \big(\textrm{cosq}^S_n(S_{\bullet})\big)_{n+1}$  soit
propre et surjectif.}  («$\textrm{cosq}_0$») de $S$, $T$ et des
différents schémas au-dessus de ceux-ci.  Le produit fibré
$X_{x,b}\times_S S^{(q)}$ (resp. ${X_T}_{x_T,b_T}\times_T T^{(q)}$) se
décompose en une somme disjointe $\coprod_{I_q} X^{(q)}_{x_i,b_j}$
(resp. $\coprod_{J_q}  {X_T^{(q)}}_{x_{Ti},b_{Tj}}$), où $I_q$
(resp. $J_q)$ est naturellement en bijection avec l'ensemble des
points $(x_i,b_j)$ de $X^{(q)}\timesf_{S^{(q)}} S^{(q)}$ (resp.
$X_T^{(q)}\timesf_{T^{(q)}} T^{(q)}$) au-dessus de $(x,b)$
(resp. $(x_T,b_T)$).  L'ensemble $I_q$ (resp. $J_q$) est donc
canoniquement en bijection avec $\pi_0(S^{(q)}_b)$ (resp.
$\pi_0(T^{(q)}_{b_T})$), où les indices désignent ici les fibres des
morphismes évidents.  (Par exemple $S^{(q)}_b$ est le produit fibré
$b\times_S S^{(q)}$, etc.)  Ainsi, les suites spectrales de descente
associées à $X_{x,b}$ et ${X_T}_{x_T,b_T}$ se réécrivent :

$$\xymatrix{
\HH^p(\coprod_{I_q}  X^{(q)}_{x_i,b_j},\mc{F}) \ar[d] \ar@2{->}[r] &
\HH^{p+q}(X_{x,b},\mc{F}) \ar[d] \\ \HH^p(\coprod_{J_q}
{X_T^{(q)}}_{x_{Ti},b_{Tj}},\mc{F}) \ar@2{->}[r] &
\HH^{p+q}({X_T}_{x_T,b_T},\mc{F})} $$ Commençons par remarquer que
pour chaque $q$, les ensembles $I_q$ et $J_q$ sont canoniquement en
bijection. C'est un cas particulier du lemme suivant :
\begin{lmm2} Soit $$\xymatrix{ \ar @{} [dr] |{\square}  T \ar[d] & T'
\ar[l] \ar[d] \\ S & S'\ar[l] }$$ un carré cartésien avec $S'\ra S$
fini. Alors, pour tout point géométrique $t$ dans $T$ d'image $s$,
$\pi_0(T'_t)\iso \pi_0(S'_s)$.
\end{lmm2} (Cela résulte par exemple du théorème de changement de base
propre ensembliste pour le morphisme fini $S'/S$ et le faisceau
$\{0,1\}$.)

Comme pour chaque $q$, les $S^{(q)}$ sont des $S'$-schémas,
l'hypothèse entraîne que  pour chaque choix de points de
$X^{(q)}\timesf_{S^{(q)}} S^{(q)}$ et  $X_{T}^{(q)}\timesf_{T^{(q)}}
T^{(q)} $ se correspondant, la flèche verticale de gauche induite sur
le facteur direct correspondant est un isomorphisme pour chaque~$p$.
\end{proof}

\section{Réduction au cas «plurinodal»}\label{reduc}

Supposons $S$ noethérien. Quitte à l'altérer $S$, on peut le supposer
également intègre. Le problème étant local en haut, on peut supposer
$f$ séparé donc $f$ compactifiable d'après un théorème de Nagata \cite{Nagata@Lutkebohmert}
($S$
est noethérien et $f$ est séparé de type fini). 
Ainsi, on se ramène au cas où $f$ est propre, quitte à prolonger $\mc{F}$ par 
zéro.

\subsection{}\label{section-descente} 
Avant d'énoncer le lemme de descente dont nous aurons besoin,
introduisons une terminologie utile ici : on dira, $f:X\ra S$, $\mc{F}$ et $r\in \ZZ$ étant donnés,
que la formation de $\cp_f(\mc{F})$ $r$-commute aux changements de base si 
les cônes des morphismes (\ref{cb}) 
n'ont pas de cohomologie en degré inférieur ou égal à $r$, pour tous les morphismes
$T\ra S$.


\begin{lmm2}\label{descente} Soient $f:X\ra S$ un morphisme propre, 
$\Lambda=\ZZ/n$ pour un entier $n$ et $r$ un entier. 
Supposons qu'il existe un $S$-hyperrecouvrement
propre tronqué à un ordre $M\geq r$,  $X_{\bullet\leq M}\ra X$, tel que la
formation des cycles proches $\cp_{X_i/S}(\Lambda)$ $(r-i)$-commute aux
changements de base, pour tout $i\in \{0,\dots,M\}$. Sous ces hypothèses,
la formation de $\cp_f(\Lambda)$ $r$-commute aux changements de base.
\end{lmm2}

\begin{proof} Commençons par remarquer que l'on peut étendre
$X_{\bullet\leq M}\ra X$  en un hyperrecouvrement propre (non tronqué)
$\varepsilon: X_{\bullet}\ra X$. Nous le faisons par commodité d'écriture.
Soient $T\ra S$ un morphisme et $(a_T,b_T)$ un
point de $T\timesf_T T$ s'envoyant sur  $(a_S,b_S)$ dans $S\timesf_S
S$. Notons $i^X_{a_S}:X_{a_S}\hra X_{S(a_S)}$ l'immersion fermée, 
$j^X_{b_S}: X_{S(b_S)}\ra X_{S(a_S)}$ le morphisme induit par localisation
étale en bas, et de la même façon les  variantes pour $X_{\bullet}$.  La
commutation aux changements de base revient à démontrer que pour tout
tel choix le cône du morphisme de changement de base (où l'on note $j_*$ pour $\R
j_*$) $${i^X_{a_S}}^* j^X_{b_S*}\Lambda_{|X_{a_T}} \ra {i^X_{a_T}}^*
j^X_{b_T*}\Lambda$$ n'a pas de cohomologie en degré inférieur ou égal à $r$. 
En effet, la fibre en $x_T$ à
gauche est isomorphe à $\cp_{f}(\Lambda)_{(x_S,a_S)}$ si $x_S$ est
l'image de $x_T$ dans $X_{a_S}$, tandis que  celle de droite est
isomorphe à $\cp_{f_T}(\Lambda)_{(x_T,b_T)}$.  Pour tout $S$-schéma
$Z$, notons $\varepsilon_Z$ l'hyperrecouvrement propre
$\varepsilon\times_S Z$ de $X\times_S Z$.  Par descente cohomologique
les morphismes $\mathrm{Id}\ra \varepsilon_* \varepsilon^*$ sont
des isomorphismes ; compte tenu également du théorème de changement de base propre
pour les $\varepsilon$, on dispose donc d'un isomorphisme canonique :
$$ {i^X_{a_S}}^* j^X_{b_S*}\Lambda\iso  {i^X_{a_S}}^*
j^X_{b_S*}\big(\varepsilon_{b_S*}\varepsilon_{b_S}^*\Lambda\big)
\isononcan \varepsilon_{a_S*} \big({i^{X_{\bullet}}_{a_S}}^*
j^{X_{\bullet}}_{b_S*}\Lambda).  $$ Le même résultat étant valable sur
$T$, on a un diagramme commutatif dont les flèches horizontales sont
des isomorphismes :

$$\xymatrix{
({i^X_{a_S}}^* j^X_{b_S*}\Lambda)_{|X_{a_T}} \ar[r] \ar[d] & \big(
\varepsilon_{a_S*} \big( {i^{X_{\bullet}}_{a_S}}^*
j^{X_{\bullet}}_{b_S*}\Lambda\big)\big)_{|X_{a_T}} \iso
\varepsilon_{a_T*} \Big( \big({i^{X_{\bullet}}_{a_S}}^*
j^{X_{\bullet}}_{b_S*}\Lambda \big)_{|X_{a_T\bullet}} \Big) \ar[d] \\
{i^{X_T}_{a_T}}^* j^{X_T}_{b_T*}\Lambda \ar[r] &  \varepsilon_{a_T*}
\big( {i^{X_{T\bullet}}_{a_T}}^* j^{X_{T\bullet}}_{b_T*}\Lambda\big)
}.$$

L'isomorphisme «en haut à droite» résulte une fois encore du  théorème
de changement de base propre et de la propreté de $\varepsilon$.

La conclusion résulte alors du lemme général suivant (appliqué au cône
des morphismes de changements de base sur les schémas simpliciaux) :
\begin{lmm2}
Soient $\varepsilon :X_{\bullet}\ra X$ un topos simplicial augmenté, $K_{\bullet}$ un complexe
de faisceaux de $\Lambda$-modules sur $X_{\bullet}$ et $r$ un entier tel que 
$\tau_{\leq r-i}({K_{\bullet}}_{|X_i})=0$ pour tout $i\in [0,r]$. Supposons 
que les topos $X_i$ aient suffisament de points.
Alors, sur $X$, le complexe de faisceaux $\tau_{\leq r} \varepsilon_* K_{\bullet}$
est trivial.
\end{lmm2}
\begin{proof}
La suite spectrale (\cite{Hodge3@Deligne},~\S~5.2.7.1)
$$E_1^{p,q}=\HH^q\varepsilon_{p*}({K_{\bullet}}_{|X_p})\Rra 
\HH^{p+q}(\varepsilon_* K_{\bullet}),$$
nous ramène à montrer que si $K'$ est un complexe de faisceaux sur $X'$
tel que $\tau_{\leq n}K'=0$ alors, pour tout morphisme de topos $e:X'\ra X$,
le tronqué $\tau_{\leq n} e_* K'$ est également nul. 
Ce dernier point est trivial. 
\end{proof}
\end{proof}

\subsubsection{}\label{descente-recurrence}

On utilisera ce lemme dans le cas où $M=r$ est strictement supérieur
au double de la dimension des fibres de $f$ ; en vertu de l'annulation de la cohomologie 
en degré $>2M$ (\ref{annulation}), la conclusion du lemme
est alors aussi forte que le théorème \ref{énoncé}.

Remarquons également qu'il est donc loisible de supposer le schéma $X$ intègre.

\subsection{}\label{dominant} Ramenons nous maintenant au cas où le
morphisme $f$ est propre et \textit{surjectif}.  On a déjà vu que l'on peut
supposer $f$ propre et $S$ intègre.  Supposons donc $f:X\ra S$ propre
et notons $F$ le fermé $f(X)\subset S$, que l'on suppose strict.
Admettons l'existence d'une modification $F'\ra F$ de $F$ telle que la
conclusion de \ref{énoncé} soit valable pour $f_{F}:X\ra F,\mc{F}$,
après changement de base à $F'\ra F$ (cela résultera de  l'étude du
cas où $f$ est dominant). Il nous reste à définir une modification de
la base $S$, partant de celle de son sous-schéma fermé $F$.

Nous allons appliquer le lemme suivant, dont l'énoncé et la
démonstration (un peu plus bas) sont dus à O.~Gabber, au morphisme
composé $g:F'\ra F\ra S$.

\begin{lmm2}\label{lemme-Ofer} Soit $g:F'\ra S$ un morphisme propre de
schémas noethériens réduits, et posons $F=g(F')$.  Il existe un
éclatement de centre nulle part dense $p:S'\ra S$ et un $F$-schéma $G$
s'envoyant par un morphisme fini et surjectif sur $p^{-1}(F)=F_{S'}$
et s'envoyant également sur $F'$.
\end{lmm2}

Soit $S'\ra S$ comme dans le lemme. Montrons que $\cp_{f_{S'}}\mc{F}$
commute aux changements de base $T\ra S'$.

Considérons le diagramme commutatif  $$\xymatrix{ & X_{F'}
\ar[dd]^{f_{F'}} \ar[dl] & & X_{G} \ar[ll] \ar[dl] \ar[dd]^{f_G} & & &
& \\ X \ar@/_1.35pc/[ddd]_f \ar[dd]^{f_{F}} & & X_{S'}
\ar@/_1.35pc/[ddd]_{f_{S'}} \ar[dd] \ar[ll] & &  X_T
\ar@/_1.35pc/[ddd]_{f_T} \ar[dd] \ar[ll] &   \\ & F' \ar[dl] & & G
\ar[ll] \ar@{->>}[dl]_{\text{fini}} & & \\ F \ar[d] & & F_{S'} \ar[d]
\ar[ll] & & F_T \ar[d] \ar[ll] & \\ S & & S' \ar[ll]_p & & T \ar[ll]
}$$ où les rectangles et les deux parallélogrammes verticaux sont
cartésiens.  La formation des cycles proches $\cp_{f_{F'}}\mc{F}$
commute aux changements de base sur $F'$.  On en déduit immédiatement
qu'il en est de même pour la formation du  complexe
$\cp_{f_{G}}\mc{F}$ relativement aux changements de base sur
$G$. D'après le lemme \ref{descente-fini}, et puisque le morphisme
$G\ra F_{S'}$ est fini et surjectif, la formation de
$\cp_{f_{F_{S'}}}\mc{F}$ commute également aux  changements de base
sur $F_{S'}$ et donc en particulier à $F_{T}\ra F_{S'}$

Il nous suffit donc d'appliquer le lemme suivant à $f_{S'}$ :

\begin{lmm2}\label{unlemme} Soient $f:X\ra S$ un morphisme propre,
d'image $F$, $n$ un entier  et $\mc{F}\in \ob \dbc(X,\ZZ/n)$.  Si la
formation des cycles proches du couple $(X\ra F,\mc{F})$ commute aux
changements de base relativement à $F$, il en est de même pour le
couple $(f:X\ra S,\mc{F})$ (relativement à $S$).
\end{lmm2}

\begin{proof} Soient $T$ un $S$-schéma et $(x_T,b_T)$ un point de
$X_{T}\timesf_T T$ d'image $(x,b)$ dans $X\timesf_S S$.  Il s'agit de
montrer que le morphisme $\RG(X_{x,b},\mc{F})\ra
\RG({X_{T}}_{x_T,b_T},\mc{F})$ est un isomorphisme.  Deux cas se
présentent : soit $b_T$ est localisé en l'ouvert complémentaire de
$F_T:=F\times_S T$, auquel cas les schémas $X_{x,b}$ et
${X_{T}}_{x_T,b_T}$ sont vides (et l'isomorphisme est évident), soit
$b_T$ est localisé en $F_T$. Dans ce dernier cas,  le morphisme
$X_{F}(x)\times_{F(a)} F(b)\ra X(x)\times_{S(a)} S(b)$ (resp. variante
sur $T$) est un isomorphisme. (En effet, $X_F=X$ et $F\times_S
S(a)=F(a)$ etc.)
\end{proof}

Il nous reste donc à démontrer le lemme~\ref{lemme-Ofer}  pour
conclure la démonstration de la réduction~\ref{dominant}.
\begin{proof}[Démonstration du lemme \ref{lemme-Ofer}] Pour chaque
sous-schéma fermé intègre $Z$ du sous-schéma fermé $F$ de $S$, il
existe un sous-schéma fermé intègre $W_Z$ de $F'$ qui est une
altération de $Z$.  (Il suffit en effet de prendre l'adhérence d'un
point fermé de la fibre de $F'\ra F$ au-dessus du point générique de
$Z$.)  D'après \ref{platification}, il existe un sous-schéma fermé
strict $R_Z$ de $Z$ tel que si l'on note $\widetilde{Z}$ et
$\widetilde{W_Z}$  les éclatés correspondants, le morphisme  induit
$\widetilde{W_Z}\ra \widetilde{Z}$ devienne fini surjectif.  Pour $Z$
variable, les ensembles localement fermés $Z-R_Z$ recouvrent $F$.  Le
schéma $F$ étant noethérien, donc compact pour la topologie
constructible, il existe un nombre fini de $Z_i,W_{Z_i}$ et $R_{Z_i}$
tels que les  $Z_i-R_{Z_i}$ recouvrent $F$.  Soient maintenant, pour
chaque tel indice $i$, $S'_i$ l'éclatement de $S$ de centre $R_{Z_i}$,
et enfin $S'$ l'éclatement de centre défini par le produit des Idéaux
$\mc{J}_{R_{Z_i}}$, \cad le «plus petit» éclatement dominant tous les
$S_i$.  Soient $(T_j)_{j}$ les composantes irréductibles de $F\times_S
S'$. Pour chaque indice $j$, choisissons un indice $i_j$ tel que le
point générique de $T_j$ s'envoie sur un point de
$Z_{i_j}-R_{i_j}$. Le schéma $T_j$ s'envoie sur le sous-schéma
$\widetilde{Z_{i_j}}$ de $S_{i_j}'$. Soit $G_j$ le produit fibré
$T_j\times_{\widetilde{Z_{i_j}}}  \widetilde{W_{i_j}}$ ; il s'envoie
par un morphisme fini et surjectif sur $T_j$ et s'envoie également sur
$F'$. L'union disjointe $G:=\coprod G_j$ répond donc à la question.

$$\xymatrix{
        &  F' \ar[d] & W_{Z_{i_j}}  \ar@{_{(}->}[l]
\ar[d]^{\mathrm{alt\acute{e}ration}} &  \widetilde{W_{Z_{i_j}}}
\ar[d]^{\mathrm{fini\,surj.}} \ar[l] & G_j \ar@{.>}[d] \ar@{.>}[l] \\
R_{i_j} \ar@{^{(}->}[r] &  F  \ar[d] & Z_{i_j}  \ar@{_{(}->}[l] &
\widetilde{Z_{i_j}} \ar[l] & T_j \ar[d] \ar@/^1.35pc/[ll]
\ar@{.>}[l]\\ &  S & & & \mathrm{\acute{E}cl}_{\cup R_i}(S) \ar[lll]
}$$

\end{proof}

\begin{rmr2} Les modifications de $F$ provenant de $S$ par restriction
ne sont pas cofinales parmi les modifications de $F$.  Il est donc
nécessaire d'avoir recourt à un morphisme intermédiaire (ici fini et
surjectif) pour palier cet inconvénient. Voici un exemple dû à János
Kollár : soit $E\hra \aff^2_{\CC}$ une courbe elliptique épointée et
$F\hra \aff^3_{\CC}=S$ le cône construit sur $E$. Soit $F'\ra F$
l'éclatement du sommet du  cône. La fibre exceptionnelle de $F'\ra F$
est non rationnellement connexe alors que les fibres de n'importe quel
éclatement de $S'\ra S$ le sont.  En particulier, $S'_F$ ne peut pas
s'envoyer par un $S$-morphisme sur $F'$.
\end{rmr2}

\subsection{}L'utilité de la réduction \ref{dominant} pour notre problème provient
du fait suivant :

\begin{prp2}\label{fibration} Soient $f:X\ra S$ un morphisme propre
\emph{surjectif} entre schémas intègres excellents, et $M$ un entier.
Il existe un schéma intègre $S'$, une altération $S'\ra S$ et un
$S'$-hyperrecouvrement propre tronqué $X'_{\bullet\leq M}\ra
X'=X\times_S S'$ tel que pour chaque $i\in [0,M]$, chaque composante
connexe du $S'$-schéma $X'_i$ soit ou bien intègre et plurinodale sur
$S'$ ou bien d'image un fermé strict de $S'$.  De plus, on peut supposer la dimension 
relative générique de $X'_0/S'$ inférieure ou égale à celle
de $X/S$. 
\end{prp2}

Rappelons qu'un morphisme est dit \emph{plurinodal} (cf.
\cite{Families@de_Jong},~5.8) s'il est le composé de morphismes projectifs et plats dont
les fibres géométriques sont des courbes connexes ayant au pire des
singularités quadratiques ordinaires.

Par commodité, nous dirons qu'un morphisme de schémas  est
\textit{presque} plurinodal s'il est somme d'un  morphisme plurinodal
et d'un morphisme  non dominant.

\begin{proof} D'après \textit{loc. cit.}~5.10, quitte à altérer $S$,
il existe une altération $X_0\ra X$ telle que $X_0/S$ soit presque
plurinodal, 
nécessairement de dimension relative égale à la dimension
de la fibre générique de $X/S$.  Le produit fibré $X_0\times_X X_0$
n'est pas nécessairement plurinodal sur $S$ ni nécessairement intègre.
Soit $Y$ le coproduit des composantes irréductibles réduites de
$X_0\times_X X_0$ ; chaque composante connexe $Y_i$ de $Y$ est donc
intègre et $Y\ra X_0\times_X X_0$ est propre et surjectif. Si $Y_i/S$
est dominant, quitte à altérer $S$ une fois de plus, on peut altérer
$Y_i$ en un $S$-schéma presque plurinodal $Y_i'$.  Dans le cas
contraire, \textit{i.e.} si $Y_i/S$ n'est pas dominant (\textit{i.e.}
surjectif), on pose $Y_i'=Y_i$. Le $S$-schéma $X_1:=\coprod Y'_i$
s'envoie par un morphisme propre et surjectif sur $X_0$ et ses
composantes connexes satisfont aux conditions de la proposition.  On
démontre alors la proposition en itérant ce procédé et de façon
répétée \cite{Hodge3@Deligne} \S~6.2.5 pour faire du schéma simplicial
\textit{strict} $(X_i)$, un vrai schéma simplicial (en particulier on
veut une flèche $X_0\ra X_1$ section de $X_1\rra X_0$ etc.).  (On
utilise implicitement le fait que la plurinodalité est une propriété
stable par changement de base.)
\end{proof}

\begin{rmr2}\label{remarque-tronquée}
Dans l'application que nous en ferons, seule l'hypothèse sur $X'_0/S'$ nous importe.
On aurait donc pu simplement considérer $\mathrm{cosq}_0^{S'}(X'_0)$, mais il n'est pas exclu
que l'énoncé ci-dessus puisse servir dans d'autres contextes.
Nous n'utiliserons donc pas le fait que les $X'_i/S'$, pour $i>0$ soient (presque) plurinodaux.
Il serait cependant intéressant et, comme l'a observé le rapporteur, non évident 
que l'on puisse supposer ces derniers de dimension relative générique inférieure à celle 
du morphisme initial $f$ ; cela suggère d'aborder cette question, qui permettrait 
d'éviter une récurrence indésirée à venir 
(cf. l'introduction du paramètre auxilliaire $r$), par des 
techniques d'hyperrésolutions cubiques
\cite{LN1335}. Rappelons que celles-ci ont l'avantage, dans le cas classique de 
la théorie de Hodge, de faire intervenir des schémas (en nombre fini) 
de dimensions de plus en plus petites (inférieures ou égales à celle du schéma de départ).
\end{rmr2}

\subsection{}\label{devissage-ultime} Reprenons dans ce paragraphe les
différentes réductions effectuées jusqu'à présent  pour voir ce qu'il
reste encore à démontrer.

\subsubsection{} D'après \ref{limite}, il suffit de démontrer le
théorème dans le cas où $S$ est un $\ZZ$-schéma de type fini ; en
particulier $S$ est noethérien, excellent, de dimension finie.\\
De plus, on peut supposer $f$ propre ; supposition qui n'est pas altérée par
les réductions faites en \S~2 et que nous ferons donc dorénavant.
Montrons alors \textit{par récurrence} sur l'entier relatif $r\geq -2$, 
et réduction au cas plurinodal et des coefficients constants, que pour tout 
$(f,\mc{F})$ comme dans \ref{énoncé}, où $\mc{F}$ est un \textit{faisceau} constructible, 
il existe une modification
$S'/S$ telle que les cônes des morphismes (\ref{cb}) n'aient pas de cohomologie en degré 
inférieur ou égal à $r$ ; cela suffit pour notre propos compte tenu de \ref{annulation}
(voir aussi \ref{descente-recurrence}). C'est évident pour $r=-2$. Soit donc $r\geq -1$
et supposons le résultat démontré pour les valeurs inférieures ou égales à $r-1$.

Pour un tel $r$, la démonstration se fait par récurrence sur la dimension $\delta_S$ de $S$ et 
la dimension relative $d_f$ de la fibre générique de $f$ (prise au-dessus de l'anneau
total des fractions si $S$ n'est pas irréductible --- ce que l'on pourrait 
supposer en vertu de \ref{descente-fini}). Amorçons la récurrence :
\begin{itemize}
\item Si $S$ est de
dimension~$0$, d'après \sga{4.5}{Th~finitude}{2.13} $f:X\ra S=s$ est universellement 
localement acyclique (relativement à tout faisceau constructible sur $X$) ; si $T$ est un
$s$-schéma strictement local, $x_T$ un point géométrique de la fibre
spéciale de $X_T\ra T$ d'image $x$ dans $X$ et enfin, $b$ un point
géométrique de $T$, on a donc :

$$
\mc{F}_x \iso \RG(X(x),\mc{F}) \lra \RG(X_T(x_T)\times_{T}
T(b),\mc{F}) \giso  \RG(X_T(x_T),\mc{F})\giso \mc{F}_{x_T}, $$ d'où la
commutation aux changements de base dans ce
cas\footnote{\label{note}Remarquons  que la définition de la locale
acyclicité donnée dans \textit{loc. cit.}  affirme que les morphismes
$\RG(X_T(x_T),\mc{F}) \ra \RG(X_T(x_T)\times_{T} b,\mc{F})$ sont des
isomorphismes. On s'intéresse ici à la cohomologie du \emph{tube} de
Milnor $X_T(x_T)\times_{T} T(b)$. Malgré tout, la démonstration donnée
dans \textit{loc. cit.} montre également l'acyclicité des tubes de
Milnor (on peut également procéder par réduction au cas lisse, par un
autre théorème de A.J. de~Jong).}.

\item Si la fibre générique de $f$ est vide, $f$ est d'image un fermé strict de $S$ de dimension
inférieure et la conclusion résulte de l'hypothèse de récurrence sur $\delta_S$ 
et de \ref{dominant} appliqué à $f(F)\subset S$.
\end{itemize}

Soient donc $f,\mc{F},n$, comme dans \ref{énoncé} avec $S$ de dimension $\delta_S$
et $f$ de dimension relative générique $d_f$, toujours supposé propre.
On peut également supposer les schémas $X$ et $S$ intègres, et $\mc{F}=\ZZ/n$ 
sans augmenter $\delta_S$ ou $d_f$. Nous laissons le soin au lecteur de vérifier
que ces réductions sont compatibles avec l'indice de récurrence
$r$ (cf. \textit{e.g.} \ref{faisceau constant} pour la dernière). Si $f$ n'est pas surjectif,
on utilise \ref{dominant} comme plus haut.

Si $f$ est surjectif, il résulte de \ref{fibration} et \ref{descente}, 
appliqués à $M=r$,
et de l'hypothèse de récurrence sur $r$ appliquée aux $X'_i/S'$ pour $i>0$
que les $r$-tronqués des morphismes de changement de base (\ref{cb}), 
sont nuls au-dessus d'une altération de $S$ convenable
pourvu qu'il en soit ainsi pour une fibration plurinodale de dimension relative égale
à $d_f$, et des coefficients constants. 

\subsubsection{}\label{rec2}Supposons donc de plus que l'on ait
démontré le théorème, au cran $r$, pour tout morphisme plurinodal $X\ra S$ de
dimension relative $d$ et pour le faisceau constant $\Lambda=\ZZ/n$.
Montrons que la conclusion est encore valable pour \emph{tout}
faisceau $\Lambda$-constructible $\mc{F}$ sur de tels $S$-schémas $X$ et donc pour tout 
\emph{complexe} de faisceaux concentré en degrés positifs. D'après les
dévissages de la seconde section, il suffit de montrer
qu'après altération, les cycles proches de $(X'\ra X \ra S,\Lambda)$
commutent aux changements de base, au sens $r$-tronqué indiqué ci-dessus, 
où $X'\ra X$ est fini, et $X'$ est
intègre.  Ici encore, deux cas se présentent : soit $X'/S$ est non
surjectif ($S$ est intègre) et l'on applique alors l'hypothèse de
récurrence sur la dimension de $S$, soit $X'/S$ est dominant et la
dimension de sa fibre générique est inférieure ou égale à
$d$. Dans ce cas, on conclut, grâce à l'hypothèse de récurrence sur $r$,
comme plus haut.

\section{Le cas plurinodal}\label{plurinodal}

Comme expliqué dans \ref{devissage-ultime}, nous allons démontrer
\ref{énoncé}, dans le cas des coefficients constants,  pour un
morphisme plurinodal, par récurrence sur sa dimension relative (notée
$d$).
On veut montrer qu'après altération les cônes de morphismes de changement de base 
\ref{cb} pour un tel morphisme n'ont pas de cohomologie en degré inférieur ou égal
à $r$ s'il en est ainsi pour $(r,d)$ plus petit (au sens lexicographique).
Remarquons que c'est seulement dans cette section que l'hypothèse faite sur la torsion
de $\Lambda$ (supposée inversible sur $S$) entre en jeu.

Le cas de la dimension relative $1$, traité en \ref{dim1}, est un
cas particulier de la proposition bien connue suivante :

\begin{prp}\label{quasi-fini} Soient $f:X\ra S$ un morphisme séparé de
type fini, $n$ un entier inversible sur $S$ et  et $\mc{F}\in \ob
\dbc(X,\ZZ/n)$ tel que le lieu de non locale acyclicité universelle de
$(f,\mc{F})$\footnote{C'est ici par définition le complémentaire du
plus grand  ouvert de $X$ sur lequel $(f,\mc{F})$  soit
\textit{localement acyclique} --- au sens de \sga{4}{xv}{1.11} --- et
le reste  après tout changement de base $T\ra S$. Une variante sans
doute plus naturelle est de considérer l'ensemble de points où $f$
satisfait ce critère ; toutefois, comme le remarque le rapporteur, il
n'est pas évident \textit{a priori} (et quoiqu'il en soit inconnu de
l'auteur) que cet ensemble soit constructible.} soit quasi-fini sur
$S$.  Alors, la formation des cycles proches commute aux changements
de base $S'\ra S$.
\end{prp}

Pour la commodité du lecteur, nous reprenons l'argument maintenant
classique de  globalisation par compactification  bête de P.~Deligne
(\textit{cf.} \sga{4.5}{Th. finitude}{}), tel qu'expliqué dans
\cite{Laumon-semi_cont},~4.1.2.

\begin{proof} On peut supposer $S$ et $S'$ strictement locaux.  Soit
$i$ (resp.~$i'$) l'inclusion de la fibre spéciale  $X_{s}\ra X$
(resp. $X'_{s'}\ra X'$). De même on définit le morphisme $j$
(resp.~$j'$), $j:X_{S(t)}\hra X$ ; notons $\cp_{s,t}$
(resp. $\cp_{s',t'}$) le foncteur $i^*\R j_*$ (resp. $i'^*\R j'_*$) ;
il s'agit d'un foncteur cycles proches «tranche par tranche».  On
dispose d'une variante évanescente évidente, notée $\ce_{s,t}$, rendue
fonctorielle si on le souhaite en travaillant dans la catégorie
dérivée filtrée adéquate.

Soit $g$ le morphisme $X'_{s'}\ra X_s$. On a une flèche de changement
de base $$g^*\cp_{s,t}\mc{F}\ra \cp_{s',t'}\mc{F}'$$  et de même pour
les cycles évanescents.  Remarquons qu'il suffit de démontrer
l'isomorphisme pour ces derniers.  À cette fin, on démontre le
résultat \textsl{a priori} plus fort suivant :
\begin{quote} Soit $x'$ un point fermé de $X'_{s'}$ d'image $x$, tel
que $x$ soit un point \textsl{isolé} de la fibre en $s$ du lieu de non
locale  acyclicité universelle de $(X,\mc{F},f)$.  Alors la flèche de
changement de base précédente est un isomorphisme  dans un voisinage
de $x'$.\end{quote}

Sous cette forme assouplie, l'énoncé est local au voisinage de $x$, si
bien que l'on peut supposer $f$ propre. La
cohomologie des  $\ce_{s,t}(\mc{F})$ (resp. $\ce_{s',t'}(\mc{F}')$)
est, au voisinage de $x$ (resp.~$x'$), concentré en ce point.  Par
propreté de $f$, si l'on applique le foncteur $\RG(X_s,-)$ au triangle
distingué

$$\trdist{\mc{F}_{\vert X_s}}{\cp_{s,t}(\mc{F}_{\vert X_{S(t)}})}{\ce_{s,t} (\mc{F}_{\vert X(s)})},$$

on obtient le triangle distingué :

$$
\trdist{\RG(X_s,\mc{F})}{\RG(X_t,\mc{F})}{\RG(X_s,\ce_{s,t}\mc{F})}.
$$

On a un morphisme évident entre ce triangle et son analogue sur $S'$ ;
à nouveau  par changement de base propre  (invariance de la
cohomologie par changement de base séparablement clos) les sommets des
arêtes horizontales sont isomorphes. Finalement, le morphisme
canonique  $$ \RG(X_s,\ce_{s,t}\mc{F})\ra
\RG(X'_{s'},\ce_{s',t'}\mc{F}) $$ est un isomorphisme.

Or, par hypothèse  $\big(\ce_{s,t}\mc{F}\big)_x$
(resp. $\big(\ce_{s',t'}\mc{F}\big)_{x'}$) est facteur direct de
$\RG(X_s,\ce_{s,t}\mc{F})$ (resp. $ \RG(X'_{s'},\ce_{s',t'}\mc{F})$),
et ces deux facteurs se correspondent.

L'isomorphisme désiré en découle.
\end{proof}

\subsection{}\label{dim1} Soit $f$ un morphisme plurinodal. Si $f$ est
de dimension relative $1$, $\mathrm{Sing}(f)$ est fini sur $S$ et
$(f,\Lambda)$ est donc universellement localement acyclique en dehors
d'un ensemble fini sur $S$ ; la proposition précédente traite ce cas
et prouve la commutation aux changements de base sans même modifier
$S$.

Supposons maintenant $f$ de dimension relative $d \geq 2$. Par
définition même,  on peut factoriser $f$ en

$$
\xymatrix{ X \ar[d]^{h\text{\ courbe\ relative}} \ar@/_1.6pc/[dd]_f \\
Y \ar[d]^g \\ S} $$ où $h$ est une courbe relative semi-stable et $g$
est plurinodal de dimension relative $d-1$.  Le complexe $
h_*\Lambda\in \ob D^{[0,2]}_c(Y,\Lambda)$ donc il existe  une modification
$S'\ra S$ telle que $\cp_{g_{S'}}(h_*\Lambda)$ commute aux changements de
base $T\ra S'$ au sens $r$-tronqué ;  cela résulte de l'hypothèse de récurrence telle que
formulée en~\ref{rec2}.  Pour simplifier les notations, supposons
$S'=S$.

Soient $T\ra S$ un morphisme, et $(a_T,b_T)$ un point de $T\timesf_T
T$ d'image $(a,b)$  dans $S\timesf_S S$ ; notons
$K^{f,T}_{a_T,b_T}(\Lambda)$ un cône du morphisme
$\cp_{a,b}(\Lambda)_{|X_{a_T}}\ra \cp_{a_T,b_T}(\Lambda)$. Comme $h$
est propre, l'image directe ${h_{a_T}}_*
\big(K^{f,T}_{a_T,b_T}(\Lambda)\big)$ (sur $Y_{a_T}$) est isomorphe au
cône $K^{g,T}_{a_T,b_T}(h_*\Lambda)$ du morphisme analogue sur $Y$ (à
coefficients  dans $h_*\Lambda$). Cette image directe n'a donc pas
de cohomologie en degré inférieur ou égal à $r$ et il en est
ainsi après tout changement  de base. On veut en
déduire qu'il en est de même de $K^{f,T}_{a_T,b_T}(\Lambda)$ ; il suffit pour cela de
savoir que le support de ce complexe est fini sur $Y_{a_T}$. Nous
allons montrer que c'est le cas, éventuellement après une modification
$S'\ra S$.

\begin{lmm}\label{unautrelemme} Sous les hypothèses précédentes, si
$\cp_g(\Lambda)$ commute au changement de base $T\ra S$, alors, pour
tout point géométrique $x_T$ de $X_{a_T}$ localisé au-dessus d'un
point  de lissité de $h$, la fibre $K^{f,T}_{a_T,b_T}(\Lambda)_{x_T}$
est nulle.
\end{lmm}

On dispose également d'une variante évidente tronquée.
Cette variante du lemme permet de conclure car, quitte à modifier $S$ une fois de
plus, on peut supposer, par hypothèse de récurrence, que
$\tau_{\leq r}K^{g,T}(\Lambda)=0$ pour tout $T/S$.  Dans ce cas, pour
tout choix de points $(a_T,b_T)$  le support de
$\tau_{\leq r}K^{f,T}_{a_T,b_T}(\Lambda)$ est inclus dans l'image inverse de
$\mathrm{Sing}(h)$ sur $X_{a_T}$, qui est fini sur $Y_{a_T}$.  Comme
$h_{a_T\, *} K^{f,T}_{a_T,b_T}(\Lambda)=K^{g,T}_{a_T,b_T}(h_*\Lambda)$, on en conclut que
$\tau_{\leq r}K^{f,T}_{a_T,b_T}(\Lambda)$ est nul ; ce qu'il fallait démontrer.

\begin{proof}[Démonstration du lemme] Soit $x_T$ comme dans l'énoncé ;
notons $x_S$ son image dans $X$ et $y_S$ l'image de $x_S$ dans $Y$.
Par lissité, le morphisme entre les localisés, $X(x_S)\ra Y(y_S)$ est
\textit{universellement} acyclique (pour les complexes constructibles
de $\Lambda$-modules).  
Par (universelle) locale acyclicité du morphisme $h$ en $x_S$, le morphisme canonique
$\RG(Y(y_S)\times_{S(a)} S(b),\Lambda)\ra \RG(X(x_S)\times_{S(a)}
S(b),\Lambda)$ est un isomorphisme ; il en est de même sur $T$.
Finalement, les fibres des cônes des flèches de changement de base
sont isomorphes  et l'on obtient :
$$K^{f,T}_{a_T,b_T}(\Lambda)_{x_T}\giso
K^{g,T}_{a_T,b_T}(\Lambda)_{y_T}=0.$$
\end{proof}

\begin{rmr}\label{remarque-adique} Le théorème \ref{énoncé} permet de
retrouver la proposition 4.2.4 de \cite{Adic_spaces@Huber} évoquée
dans l'introduction.  La raison principale en est qu'une modification
d'un schéma valuatif $S$ possède une section (cf. \textit{e.g.}
\ega{ii}{7.3.1}) si bien que  l'on a commutation aux changements de
base sur $S$.  Notons à ce propos que la démonstration de Roland
Huber, comme celle de P.~Deligne, n'entraîne pas \emph{a priori} la
commutation aux changements de base $T\ra S$ généraux mais seulement
ceux où $T$ est un schéma valuatif (dominant $S$).  (Leurs résultats
démontrent cependant l'injectivité des morphismes de changements de
base pour chaque groupe de cohomologie $\HH^i \cp$.)
\end{rmr}

\begin{rmr} Comme le remarque O.~Gabber, la conclusion du théorème
\ref{énoncé} vaut également si $f$ est seulement supposé de type fini
mais que l'espace topologique sous-jacent à $S$ est noethérien. Si $S$
est affine --- cas auquel on s'est ramené en \ref{limite} ---, cela
résulte du fait que $X_{\red}$ est aussi le schéma réduit associé à un
$S$-schéma de présentation finie.  (Rappelons que l'espace topologique
sous-jacent à $S=\SP(A)$ est  noethérien si et seulement si le radical
de tout idéal $I$ de $A$ est égal au radical d'un idéal de type fini ;
on utilise alors le fait démontré dans  \cite{noetherian@Ohm} que si
$X$ est affine de type fini sur un tel $S$, l'espace topologique
sous-jacent à $X$ est également noethérien.)

\end{rmr}

\nocite{ega,sga4xv}

Terminons cette partie par un complément.

\section{Calcul des fibres}

Étant donné un morphisme $f:X\ra S$ et un faisceau $\mc{F}$ sur $X$,
on considère aussi  classiquement (cf.
\sga{4.5}{Th. finitude}{2.11} et note en bas de page
(\textit{\ref{note}}))  la cohomologie des \textit{fibres} de Milnor
\mbox{$X(x)\times_{S(s)} t$} à  valeur dans $\mc{F}$ en plus de celle
des \textit{tubes}  \mbox{$X(x)\times_{S(s)} S(t)$}, considérée dans
la première partie.  Il résulte immédiatement du théorème \ref{énoncé}
qu'après modification de la base, elles sont canoniquement
isomorphes. En fait on peut démontrer un résultat plus précis :

\begin{thm}\label{tubes-fibres}\label{cohopropre} Soient $f:X\ra S$ et
$\mc{F}$ comme en \ref{énoncé}. Pour tout $N\in \NN$,  il existe une
modification $S'_N\ra S$ de $S$ telle que pour tout $S'_N$-schéma $T$,
et tout point géométrique $x_T$ de $X_{T}$, le couple
$({f_{T}}_{x_T},\mc{F}_{|X_T})$  soit \emph{cohomologiquement propre}
en degré $\leq N$.

En particulier, après modification de la base,  la cohomologie des
tubes de Milnor s'envoie isomorphiquement sur la cohomologie des
fibres de Milnor : pour tout point  $(x_T,b_T)$ de $X_T\timesf_{T} T$,
d'image $(a_T,b_T)$ dans $T\timesf_T T$, le morphisme $$
\RG(X_T(x_T)\times_{T(a_T)} T(b_T),\mc{F})\ra
\RG(X_T(x_T)\times_{T(a_T)} b_T,\mc{F}) $$ est un isomorphisme.
\end{thm}

Rappelons que la propreté cohomologique en degré $\leq N$ signifie que
la formation de l'image directe tronquée $\tau_{\leq N}\big(\R
{f_T}_{x_T*}\mc{F}\big)$ commute aux changements de base $Z\ra
T(s_T)$, où $s_T$ est l'image du point géométrique $x_T$ par $f_T$.
On conjecture (cf. remarque ci-dessous) que la variante non tronquée
du théorème est également vraie.
Les arguments étant les mêmes que ceux présentés dans la partie
précédente, nous nous contentons d'une esquisse de démonstration.

\begin{proof} Le second énoncé (tube versus fibre) résulte du premier
(ou de \ref{énoncé} comme  annoncé plus haut) car on sait que ni
$\RG(X_T(x_T)\times_{T(a_T)} T(b_T),\mc{F})$, ni
\mbox{$\RG(X_T(x_T)\times_{T(a_T)} b_T,\mc{F})$} n'ont de cohomologie
en degré supérieur à deux fois un majorant strict de la dimension des
fibres de $f$. (On utilise \ref{annulation} et, en passant à la
limite, la dimension cohomologique des schémas de type fini sur
$\SP\,\kappa(b_T)$.)

Soient $s$ un point géométrique de $S$ et $j_Z: Z\ra S(s)$ un
morphisme  entre schémas strictement locaux. Notons $z$ le point fermé
de $Z$, $S(t)$ le localisé strict de $S(s)$ en $t:z\ra S$ et $j_{S,z}$
le morphisme de localisation  $S(t)\ra S(s)$. Soit $i_s$ (resp.~$i_z$)
l'immersion fermée $s\hra S(s)$ (resp. $z\hra Z$).  Notons avec un $X$
en exposant, les morphismes obtenus par produit fibré avec $X\ra S$.
Il s'agit donc de montrer que pour chaque entier $N$, quitte à
modifier $S$,  pour tout choix de $s$ et $Z$ comme précédemment, le
morphisme d'adjonction 
\begin{equation}\label{eqpropcoh} \big({i^X_s}^*
j^X_{S,z\,*}\mc{F}\big)_{|X_z}\ra {i^X_z}^* j^X_{Z\,*} \mc{F}
\end{equation} (entre faisceaux sur $X_z$) est un isomorphisme en
degré $\leq N$ et qu'il en est encore ainsi pour tous les $S$-schémas
$T$.  En effet, la fibre en un point géométrique $x$ de $X_z$ du
morphisme \ref{eqpropcoh} s'identifie au morphisme de changement de
base pour le carré cartésien $$\xymatrix{ \ar @{} [dr] |{\square}
X_{S(s)}(x) \ar[d]^{f_x} & X_{S(s)}(x)\times_{S(s)} Z \ar[l]_{j_Z^X}
\ar[d]^{f_x\times_{S(s)} Z} \\ S(s) & Z \ar[l]_{j_Z}  \\ t & z
\ar@{|->}[l] } $$ \cad :

$$
\RG(X_{S(s)}(x)\times_{S(s)} S(t),\mc{F})=(f_{x*}\mc{F})_t \ra
\RG(X_{S(s)}(x)\times_{S(s)} Z,\mc{F})=\big(
(f_{x}\times_{S(s)}Z)_*\mc{F}\big)_z .$$

Notons ${K'}^{f,Z}(\mc{F})$ un cône de \ref{eqpropcoh}.  Reprenant les
dévissages de I.\ref{réductions}, on vérifie sans peine que l'on peut
supposer $\mc{F}=\Lambda$ et $f$  plurinodal\footnote{Par exemple, on
démontre l'analogue du lemme \ref{unlemme} en remarquant que si l'on
factorise un morphisme propre $f:X\ra S$ ($S$ strictement local) en
$X\sr{g}{\ra} F \sr{i}{\hra} S$, alors pour tout point géométrique $x$
de la fibre spéciale, on dispose d'isomorphismes  $\R
f_{x*}\mc{F}\isononcan i_* \R g_{x*}\mc{F}$ et $\R
(f_{x}\times_{S}Z)_*\mc{F}\isononcan (Z\times_S F \hra Z)_*  \R
(g_{x}\times_{S}Z)_*\mc{F}$.}. Pour alléger l'exposition, nous omettons
la récurrence sur le degré de troncation $r$, identique à celle effectuée pour
démontrer le théorème de changement de base. Le lecteur la rétablira aisément
de lui-même.

Pour pouvoir procéder par récurrence (suivant \ref{rec2}) dans le cas
plurinodal, comme expliqué dans la partie précédente, il nous suffit
de vérifier trois propriétés de l'obstruction $K'$.

\begin{itemize}
\item L'obstruction $K^{'\,Z}$ commute aux images directes propres.
(C'est vrai d'après le théorème de changement de base propre.)

\item L'obstruction $K^{'\,Z}$ est nulle sous les hypothèses de
\ref{quasi-fini}, \cad lorsque le lieu de non locale acyclicité est
quasi-fini sur $S$.  (En effet, en un point de lissité $x$ de $f$, les
morphismes $\Lambda\ra \RG(Z,\Lambda)\ra \RG(X(x)\times_{S(s)}
Z,\Lambda)$ sont des isomorphismes, ce qui nous permet d'utiliser le
même argument par compactification.)

\item Avec les notations et les hypothèses de \ref{unautrelemme}, on a
des isomorphismes ${K'}^{f,Z}(\Lambda)_x\giso {K'}^{g,Z}(\Lambda)_y
=0$. (En un point de lissité le morphisme induit entre les localisés
est \emph{universellement} localement acyclique.)
\end{itemize}
\end{proof}

\begin{rmrs} On devrait disposer, \emph{après modification de la
base},  d'un analogue de la proposition \ref{annulation} pour les complexes
$\RG(X(x)\times_{S(s)} Z,\mc{F})$.  Alternativement, il est
vraisemblable, comme le suggère O.~Gabber, que l'énoncé de propreté
cohomologique se ramène au cas où $Z$ est le spectre d'un anneau de
valuation dont le corps des fractions est algébriquement clos ; sous
cette hypothèse,  O.~Gabber sait démontrer le résultat d'annulation
requis. Cela  permettrait donc d'éviter le recourt à des énoncés
tronqués.\\  Enfin, remarquons qu'il aurait été possible de donner une
démonstration  uniforme du théorème \ref{énoncé} et du théorème
précédent mais l'axiomatisation qu'elle requiert ne semble pas de
nature à éclaircir l'argument.
\end{rmrs}

\subsection{}\label{ubiquité-traits} Soient $S$ un schéma noethérien,
$a$ un point géométrique  de $S$ et $b$ une générisation géométrique
de $a$.  Il existe un trait $T_{a,b}$ strictement local et  un
morphisme $T_{a,b}\ra S$ tel que  $a$ (resp.~$b$) soit dominé par le
point fermé (resp.  l'image du point générique géométrique) de
$T_{a,b}$ (cf. \ega{ii}{7.1.7}).  Ainsi, on peut compléter le
théorème précédent, de la façon suivante : sous les hypothèses de
\ref{énoncé}, quitte à modifier la base $S$, chaque fibre des cycles
proches $\cp_{f}(\mc{F})$ est isomorphe à la fibre d'un complexe de
cycles proches usuels.  Plus précisément, étant donné $(x,b)$, d'image
$(a,b)$ dans $S\timesf_S S$,  il existe un trait $T_{a,b}\ra S$ et un
point géométrique $\tilde{x}$ de $X_{T_{a,b}}$ au-dessus de $x$ tel
que l'on ait des isomorphismes : $$ \cp_{f}(\mc{F})_{(x,b)}\isononcan
\RG(X(x)\times_{S(a)} S(b),\mc{F}) \iso
\RG(X_{T_{a,b}}(\tilde{x})\times_{T_{a,b}} b,\mc{F})\isononcan
\cp_{f_{T_{a,b}}}^{\mathrm{usuel}}(\mc{F})_{\tilde{x}}, $$ où
$\cp_{f_{T_{a,b}}}^{\mathrm{usuel}}(\mc{F})$ le complexe des cycles
proches défini par A.~Grothendieck dans \sga{7}{I}{}.

Cela découle de \ref{énoncé}.

\part{Constructibilité}

\section{Énoncés}

Il résulte de \ref{ubiquité-traits} et du théorème de constructibilité
de P.~Deligne (\sga{4.5}{Th. finitude}{3.2}), qu'après modification de
la base, les fibres des cycles proches sont finies, et même qu'ils sont constructibles sur chaque fibre,
c'est-à-dire une fois restreints à $X_s\timesf_{S(s)} t$, où $(s,t)$
est un point de $S\timesf_S S$.  On a en fait un résultat plus précis :

\begin{thm}\label{constructibilite} Soient $f:X\ra S$ et $\mc{F}$
comme en \ref{énoncé}.  \emph{Quitte à modifier $S$}, on peut supposer
que pour tout entier naturel $i$, les faisceaux $\HH^i\big( \cp_{f}
\mc{F}\big)$ sont constructibles sur $X \timesf_{S} {S}$, \cad qu'il
existe des partitions finies de $X$ et $S$ en des parties
constructibles  localement fermées, $X=\bigcup X_{\alpha}$ et
$S=\bigcup S_{\beta}$, telles que chaque $\HH^i\big( \cp_{f}
\mc{F}\big)$ soit localement constant fini sur les sous-topos
$X_{\alpha}\timesf_{S} S_{\beta}$ de $X\timesf_S S$.
\end{thm}

\begin{crl}\label{constructibilite2} Soient $f:X\ra S$ et $\mc{F}$
comme en \ref{énoncé}.  Il existe une modification $S'\ra S$ telle que
pour tout $S'$-schéma $T$, et tout point géométrique $x_T$ de $X_T$,
d'image $a_T$ dans $T$, le \emph{complexe} ${f_T}_{x_T\,*} \mc{F}$ soit
constructible  sur $T(a_T)$.
\end{crl}

\begin{rmr}
Il résulte de \sga{4}{xvii}{5.2.11} et \ref{annulation} que le foncteur $\cp_f$ préserve la
propriété d'être de Tor-dimension finie. Joint à \ref{constructibilite},
cela entraîne, pour un $f$ comme dans \textit{loc. cit.} et un nombre
premier $\ell$ inversible, qu'après modification
de la base, les complexes $\cp_f(\ZZ/\ell^{n+1})$ sont dans $\mathsf{D}_{ctf}(X\timesf_S S,\ZZ/\ell^{n+1})$
pour tout $n\geq 0$.
\end{rmr}

Comme précédemment, on se ramène au cas où $S$ est noethérien de
dimension finie, excellent et intègre. On procède par récurrence sur
$\dim(X/S)$ et $\dim(S)$.

La démonstration se fait en deux étapes : passer du cas plurinodal au
cas général puis démontrer la constructibilité dans le cas plurinodal
(et des coefficients constants).
Avant cela, commençons par quelques propriétés de ces topos et des
faisceaux constructibles sur iceux, en suivant les suggestions du
rapporteur\footnote{L'auteur reste bien sûr seul responsable des
erreurs qui pourraient se trouver dans se texte.}.

\section{Conditions de finitude dans les topos $X\timesf_Y Z$}

\subsection{Cohérence ; rappels}\label{coherence}Remarquons tout
d'abord que  la définition donnée dans \ref{constructibilite} garde un
sens pour les faisceaux d'ensembles ou de $\Lambda$-modules sur des
topos de la forme $E=X\timesf_Y Z$, où $X,Y$ et $Z$ sont des schémas
cohérents. Comme dans \cite{Vanishing@Laumon}, on définit un site pour
le topos $E$, dont les objets sont les triplets $(U_X/V_Y\backslash
W_Z)$ où  $U_X,V_Y,W_Z$ sont étales séparés de type fini sur $X,Y,Z$.

Un tel objet définit naturellement un faisceau d'ensembles, qui est
constructible.  En effet, si $X'$ est un sous-schéma de $X$
(resp.~$Z'$ un sous-schéma de $Z$) tel que la restriction de $U_X$ à
$X'$ (resp. de $W_Z$ à $Z'$, resp. l'image inverse de $V_Y$ sur $Z'$)
soit \textit{finie} étale,  la restriction du faisceau représenté par
$U_X/V_Y\backslash W_Z$ au sous-topos $X'\timesf_Y Z'$ est localement
fini constant.

On peut également étendre au présent cadre la proposition
\sga{4}{ix}{2.7}, affirmant que de tels faisceaux sont  des
générateurs de la catégorie des faisceaux constructibles (avec une
variante pour les $\Lambda$-modules). Lors de la démonstration du
théorème \ref{constructibilite}, une caractérisation plus maniable des
faisceaux constructibles sera utile ; dans le cas de topos
raisonnables, il s'agit des objets de \textit{présentation finie}
(\textit{i.e.} les objets $\mathrm{PF}$ tels que $\Hom(\mathrm{PF},-)$ commute
aux colimites filtrantes), voire également des objets
\textit{noethériens}.  Ces questions, dans le cas des topos comme $E$,
font l'objet du reste du paragraphe.

Commençons par remarquer le topos $E$ est \textit{cohérent}
\sga{4}{vi}{2.3}, car le site de définition introduit ci-dessus
satisfait aux conditions de la proposition 2.1 de \textit{loc. cit.}
et possède un objet final.

De plus, un faisceau constructible d'ensembles $\mc{F}$ sur $E$  est
nécessairement  un objet cohérent de $E$ (au sens de
\sga{4}{vi}{1.13}, \cad ici encore quasi-compact et quasi-séparé
(cf. \textit{loc. cit.})).  Cela résulte d'une part du fait qu'un
objet localement fini constant d'un topos cohérent est cohérent et
d'autre part du fait que si l'on a un morphisme cohérent $p:E'\ra E$
de topos cohérents et si $p$ est conservatif, alors pour tout objet
$\mc{F}$ de $E$, $\mc{F}$ est cohérent si et seulement si $p^*\mc{F}$
l'est.  (On applique le deuxième point au topos $E'$ somme disjointe
(finie) des sous-topos de $E$ sur lesquels le faisceau $\mc{F}$ est
localement fini constant et $p$ le morphisme canonique correspondant.)
On en déduit que tout objet de $E$ est colimite filtrante d'objets
cohérents, si bien que le topos $E$ est \textit{parfait}, au sens de
\sga{4}{vi}{2.9.1}. Les notions topossiques de constructibilité,
cohérence et de présentation finie (\textit{loc. cit.}~1.9.3, 1.24)
coïncident ici.

Il reste à vérifier que les objets constructibles sont également
noethériens.

\subsection{Localisation}Rappelons qu'en vertu d'un théorème de
P.~Deligne, nos topos, qui sont localement cohérents, ont suffisamment
de points (\textit{loc. cit.},~9.0). Étant donné un point $p$ d'un
topos $T$, on peut former le topos \textit{localisé en} $p$,
la $2$-limite $\lim_{U\ni p} T_{/U}$, où $U$ parcourt les voisinages
ouverts de $p$.  C'est un topos \textit{local} au sens de
\sga{4}{vi}{8.4} : le foncteur section globale est un foncteur fibre,
de \textit{centre} $p$, noté $T_p$ («localisation de $T$ en $p$»).
Remarquons que le topos $T\timesf_T T$ joue le rôle de localisation
universelle : on a un diagramme $2$-cartésien $$ \xymatrix{\ar @{}
[dr] |{\square} T_p \ar[d] \ar[r] & T\timesf_T T
\ar[d]^{\mathrm{pr}_1} \\ p \ar[r] & T} $$

\subsection{Quelques lemmes}

\begin{lmm2} Soient $X\ra T \la T' \la S$ des morphismes de schémas
avec $T'\ra T$ entier. Posons $X'=X\times_T T'$.  Le morphisme
$p:E'=X'\timesf_{T'} S\ra X\timesf_T S=E$ est une équivalence de topos.
\end{lmm2}

\begin{proof} Il suffit de montrer que les foncteurs $p_*$ et $p^*$
sont quasi-inverses l'un de l'autre. Compte tenu du fait que les topos
considérés ont suffisamment de points, cela résultera du fait que les
flèches d'adjonctions sont des isomorphismes sur les fibres. Or, tout
point $(x,s)$ de $E$ se relève uniquement en un point $(x',s)$ de $E'$
et si $E$ est local centré en $(x,s)$, $E'$ est également local et
centré en $(x',s)$. Le résultat en découle.
\end{proof}

L'ensemble des classes d'isomorphismes de points d'un topos $E$ comme
plus haut peut-être ordonné par la relation d'ordre : $p\leq q$ si et
seulement si il existe un morphisme $q\ra p$ (i.e. $p$ est une
spécialisation de $q$). (On remarquera que tout endomorphisme d'un
point est un automorphisme donc cette relation est bien antisymétrique.)

\begin{lmm2} Soit $f:X\ra S$ un morphisme de schémas noethériens avec
$f$ de type fini ou bien $S$ excellent.  Alors :
\begin{enumerate}
\item L'ensemble des classes d'isomorphismes de points de $E=X\timesf_S S$ a un nombre
fini d'élément maximaux,
\item Toute classe d'isomorphisme est majorée par un élément maximal.
\end{enumerate}
\end{lmm2}

\begin{proof} Tout point de $E$ est spécialisation d'un point
au-dessus d'un point générique d'une composante irréductible de
$S$. Le schéma $S$ étant noethérien, il suffit donc de considérer le
cas de $E'=X\timesf_S s$ où $s$ est un schéma ponctuel de $S$.  Soit
$S'$ la normalisation de l'adhérence de $s$ dans $S$ et $X'$ le
produit fibré $X\times_S S'$.  D'après le lemme précédent, $E'$
s'identifie à $X'\timesf_{S'} s$. Chaque point de $X'$ est  dominé par
une unique classe d'isomorphisme de points de $E'$ ; en effet, si $t$
est un point de $S'$, $s'$ et $t'$ des points topossiques au-dessus de
$s$ et $t$, le groupe $\Aut(s')$ agit transitivement sur
$\Hom(s',t')$.  La relation d'ordre sur les points de $E'$ se déduit
de celle des générisations au sens classique (Zariski) sur $X'$. Il
suffit donc de montrer que le schéma $X'$ n'a qu'un nombre fini de
composantes irréductibles. Si $f$ est de type fini, on utilise le
fait, démontré dans \cite{NoetherianSpectrum@Heinzer}, que la
normalisation d'un anneau intègre noethérien a un spectre
noethérien. Dans le cas où $S$  est excellent on conclut plus
simplement par le fait que la normalisation est finie sur la base.

Comme il résulte de la démonstration, le second énoncé est même
valable sans hypothèses.
\end{proof}

\begin{lmm2} Sous les hypothèses du lemme précédent, on a :
\begin{enumerate}
\item Tout sous-faisceau d'un faisceau constructible est constructible,
\item Tout faisceau constructible sur $E$ est noethérien.
\end{enumerate}
\end{lmm2}

\begin{proof} 
Commençons pas démontrer qu'une suite croissante de sous-faisceaux
constructibles $(\mc{F}_n)$ d'un faisceau constructible $\mc{F}$ sur $E$ est stationnaire.
Par récurrence noethérienne, on peut supposer que la
conclusion est valide pour $X'\ra S'$ où $X'$ et $S'$ sont des
sous-schémas fermés de $X,S$, avec $(X',S')$ différent de $(X,S)$.
Il suffit de montrer qu'elle est stationnaire au-dessus d'un sous-topos
$E'=U_X\timesf_S V_S$ pour des ouverts denses $U_X\hra X$ et $V_S\hra S$.
Il existe un entier $N$ tel que pour chaque point maximal $p$, les suites
${\mc{F}_n}_p$ sont constantes pour $n\geq N$. Si l'on prend $U_X$ et $V_S$
tels que $\mc{F}_N$ et $\mc{F}$ soient localement constants sur $E'$, cela entraîne
les isomorphismes $\mc{F}_n=\mc{F}_N$ sur $E$ pour $n\geq N$.
Prouvons 1). Soit $\mc{G}$ est un sous-faisceau d'un faisceau constructible $\mc{F}$.
D'une part il est colimite filtrante de faisceaux constructibles, et d'autre part,
l'image d'un morphisme de faisceaux constructibles est constructible ; le faisceau
$\mc{G}$ est donc également réunion d'un ensemble filtrant croissant de 
sous-faisceaux \textit{constructibles} de $\mc{F}$. D'après ce qui précède 
cet ensemble a un élément maximal et $\mc{G}$ est donc constructible.
Le second point en résulte immédiatement (par exemple en l'appliquant à la «limite»
d'une hypothétique suite strictement croissante).

\end{proof}

\section{Démonstration du théorème \ref{constructibilite}}\label{dem-constr}

\subsection{Réduction au cas plurinodal}

Commençons par un lemme général, qui trouverait certainement sa place
plus haut.

\begin{lmm2}\label{calcfibre} Soient $h:X\ra Y$ un morphisme propre
entre $S$-schémas et $\fl{\mc{F}}$ un faisceau de torsion sur
$X\timesf_S S$. Si $(y,t)$ est un point de $Y\timesf_S S$, il existe
un morphisme canonique de topos $\varphi_{y,t}:X_y\ra X\timesf_S S$,
«$x\mapsto (x,t)$», tel que la fibre de $\R \fl{h}_* \fl{\mc{F}}$ en
$(y,t)$ soit canoniquement isomorphe à
$\RG(X_y,\varphi_{y,t}^*\fl{\mc{F}})$.
\end{lmm2}

\begin{proof} Soit $(y,t)$ comme dans l'énoncé et notons $s$ l'image
de $y$ dans $S$. Quitte à remplacer $S$ par $S(s)$ on peut supposer
$S$ strictement local. De même, on peut supposer $Y(y)=Y$.  La fibre
de $\R \fl{h}_* \fl{\mc{F}}$ en $(y,t)$ s'identifie naturellement à
$\RG(X\timesf_S S(t),\mc{F})$. Notons $\pi_{t}$ le morphisme
$X\timesf_S S(t)\ra X$ ; le théorème de changement de base propre
permet de réécrire $\RG(X\timesf_S
S(t),\mc{F})\isononcan\RG(X,\R{\pi_t}_*\mc{F})$ comme
$$\RG(X_y,i_y^*\R{\pi_t}_*\mc{F}),$$ où $i_y$ est l'immersion fermée
de la fibre spéciale $X_y\hra X$.  Remarquons maintenant qu'il existe
un morphisme de topos $\varphi_{y,t}$  et un diagramme commutatif : $$
\xymatrix{ X\timesf_S S(t) \ar[d]_{\pi_t} & \\ X & X_y \ar[l]_{i_y}
\ar@_{.>}[ul]_{\varphi_{y,t}} } $$ au-dessus du morphisme $s=\Ens\ra
S\timesf_S S(t)$, $s\mapsto (s,t)$.  Le point clé indiqué à l'auteur
par O.~Gabber,  dont le lemme résulte immédiatement, est que le
morphisme canonique ${i_y}^*\R {\pi_t}_*\ra \varphi_{y,t}^*$ est un
isomorphisme. En passant aux fibres, il nous suffit de montrer que si
$x$ est un point géométrique de $X_y$, le morphisme canonique
$\RG(X(x)\timesf_S S(t),\mc{F})\ra \mc{F}_{(x,t)}$ est un
isomorphisme.  Cela résulte du fait que le topos $X(x)\timesf_S S(t)$
est un topos local, de centre $(x,t)$.
\end{proof}

Ce lemme entraîne en particulier que la formation des images directes
par un morphisme propre «fléché» commute aux changements de base
$S'\ra S$ et que si $\mc{F}$ est un faisceau sur $X$, on 
a un isomorphisme $\fl{h}_*p^*\mc{F}
\iso p^*\fl{h}_*\mc{F}$ (où les $p$ sont les premières projections
des topos fléchés vers les topos usuels). 

\begin{prp2}\label{constrimagedirecte} Soient $S$ un schéma
noethérien, $f:X\ra Y$ un morphisme propre entre $S$-schémas de type
fini et $\fl{\mc{F}}$ un faisceau de $\Lambda$-modules constructible
sur $X\timesf_S S$. Alors, l'image directe $\R \fl{f}_* \fl{\mc{F}}$
est constructible sur $Y\timesf_S S$.
\end{prp2}

Remarquons qu'il n'est pas nécessaire, tout comme dans le cas
classique (auquel nous nous ramenons), que la torsion de
$\Lambda=\ZZ/n\ZZ$ soit d'ordre inversible sur $S$.

\begin{proof} Pour tout objet $\mc{U}$ du site définissant
$X\timesf_S S$ introduit en \ref{coherence}, notons
$\Lambda_{\mc{U}/X\timesf_S S}$ (ou plus simplement
$\Lambda_{\mc{U}}$) le faisceau associé à $\Lambda[\Hom(-,\mc{U})]$
sur $X\timesf_S S$. Le faisceau $\Lambda$-constructible, donc
noethérien, $\fl{\mc{F}}$ est isomorphe à un quotient des sommes
directes finies de tels faisceaux (cf.
\textit{e.g.} \sga{4}{ix}{2.7}).  Prenant une résolution à gauche de
$\fl{\mc{F}}$ par somme directes de faisceaux de le forme
$\Lambda_{\mc{U}}$ et considérant la première suite spectrale, on se
ramène à démontre la proposition dans le cas particulier d'un faisceau
$\fl{\mc{F}}$ de la forme $\Lambda_{\mc{U}}$. (On utilise
implicitement le fait que la catégorie des faisceaux constructibles
est stable par (co)noyaux et extensions.)

Remarquons qu'il suffit de démontrer le théorème après un changement
de base fini et surjectif  $S'\ra S$ (cf. lemme
\ref{innocuite}). Or, si l'ouvert $\mc{U}$ de  $X\timesf_S S$ est de
la forme $U_X / V_S \backslash U_S$, il existe un tel morphisme $S'\ra
S$ tel que l'image inverse de $\mc{U}$ sur $X_{S'}\timesf_{S'} S'$
soit un coproduit d'ouverts de la forme $U_{X'}/ V_{S'} \backslash
U_{S'}$ où maintenant $V_{S'}$ et $U_{S'}$ sont des ouverts de Zariski
de $S'$. Plus précisément on a :

\begin{lmm2} Soit $f:U\ra S$ un morphisme séparé étale entre schémas
cohérents. Il existe un morphisme fini surjectif $S'\ra S$ tel que
l'image inverse de $U$ sur $S'$ soit un coproduit d'ouverts de Zariski.
\end{lmm2}

\begin{proof} Par le Main Theorem de Zariski, on peut factoriser $f$
en une immersion ouverte suivie d'un morphisme fini  $f:U\hra T\ra
S$. Quitte à rajouter à $T$ le complémentaire de l'image de $f$, on
peut supposer $T\ra S$ également surjectif. Soit $W$ le sous-schéma
ouvert de $S$ au-dessus duquel le cardinal des fibres géométriques de
$f$ est maximal, égal à $m$. Le morphisme $f$ est fini sur $W$ et si
$m>0$, on peut supposer $U=T$ sur $W$. L'image inverse de $U$ sur $T$
se décompose ; on conclut alors en itérant ce procédé et avec $m$ de
plus en plus petit.
\end{proof}

Pour simplifier les notations supposons dorénavant $S'=S$.  Pour des
ouverts $\mc{U}$ comme ci-dessus, il est évident que la
constructibilité de $\R \fl{f}_* \big(\Lambda_{\mc{U}}\big)$ résulte
de la constructibilité de $\R \fl{f}_* \Lambda_{U_X/S\backslash S}$
sur $Y\timesf_S S$.  Si l'on note $\mathrm{pr}_1$ la projection
$X\timesf_S S$, on a un isomorphisme canonique
$$\Lambda_{(U_X/S\backslash S) \, / \, (X\timesf_S
S)}\isononcan\mathrm{pr}_1^*\Lambda_{U_X \, / \, X}$$ où
$U_X/S\backslash S$ est l'image inverse sur $X\timesf_S S$ de l'ouvert
$U_X$ de $X$.

Le faisceau $\Lambda_{U_X/X}$ étant constructible sur $X$ il est
isomorphe à sous-faisceau d'un produit de faisceaux $\pi_*C$ où
$\pi:X'\ra X$ est un morphisme fini et $C$ un $\Lambda$-module
constant constructible.  Il nous reste donc à montrer que $\R \fl{f}_*
(\mathrm{pr}_1^*\pi_*\Lambda)$ est constructible sur $Y\timesf_S
S$. Ce dernier est isomorphe à $\R \fl{f'}_*\Lambda$, où $f'$ est le
composé $X'\ra X\ra Y$, lui-même isomorphe à ${\mathrm{pr}'_1}^* \R
f'_* \Lambda$. D'après le théorème  de constructibilité classique, $\R
f'_* \Lambda$ est constructible sur $Y$, ce qui achève la
démonstration de la proposition.
\end{proof}

\begin{lmm2}\label{innocuite} Soit $(f,\fl{\mc{F}})$ comme dans la
proposition précédente et supposons qu'il existe un morphisme  fini
surjectif $S'\ra S$  tel que $\R \fl{f_{S'}}_* \fl{\mc{F}}_{S'}$ soit
constructible sur $Y_{S'}\timesf_{S'} S'$. Alors, $\R \fl{f}_*
\fl{\mc{F}}$ est constructible sur $Y\timesf_S S$.
\end{lmm2}

Compte tenu de la commutation aux changements de base, ce lemme
résulte à son tour du lemme :

\begin{lmm2}\label{innocuite2} Soient $S'\ra S$ un morphisme propre
surjectif, $X$ un $S$-schéma, et $X'$ son image inverse sur $S'$.
Considérons $p:E'=X'\timesf_{S'} S'\ra X\timesf_{S} S=E$. Alors, si
$\fl{\mc{F}}$ est un faisceau sur $E$ tel que  $p^*\fl{\mc{F}}$ soit
constructible sur $E'$, le faisceau $\fl{\mc{F}}$ est également
constructible.
\end{lmm2}

\begin{proof} On remarque que tout point de $E$ peut se relever en un
point de $E'$ si bien que le morphisme $p$ est conservatif. Il est
aussi cohérent, comme on le voit sur les sites de définition. On
utilise alors le fait, énoncé dans le dernier paragraphe de
\ref{coherence}, que la cohérence (équivalente à la constructibilité)
se teste après image inverse par un morphisme cohérent conservatif.
\end{proof}

La proposition précédente étant établie, il résulte des dévissages
utilisés pour établir la commutation aux changements de base
(cf. plus particulièrement \ref{devissage-ultime}) qu'il
suffit de démontrer le théorème dans le cas des coefficients constants
et d'une fibration plurinodale. Ici encore, compte tenu de \ref{remarque-tronquée},
il nous faut procéder par récurrence sur l'entier $r$ tel que 
les $\tau_{\leq r}\cp_{X/S}(\mc{F})$ ($X/S$ et $\mc{F}$ variables) deviennent constructibles
au-dessus d'une altération de $S$.
Plus précisément, si $\varepsilon:X_{\bullet}\ra X$ est un hyperrecouvrement propre de $X$
comme dans \ref{fibration} (pour $M=r$ et étendu arbitrairement
en plus grand degré) et $\fl{\varepsilon}$ est le morphisme de
topos induit par $\varepsilon$ entre le topos simplicial associé aux
$X_i\timesf_S S$ (noté $X_{\bullet}\timesf_S S$) et $X\timesf_S S$, on a un isomorphisme
canonique entre $\cp_{X/S}(\Lambda)$ et $
\fl{\varepsilon}_*\cp_{X_\bullet/S}(\Lambda)$.  Les constituants de
$\cp_{X_{\bullet}/S}(\Lambda)$ ne sont autres que les
$\cp_{X_i/S}(\Lambda)$ ($i\geq 0$).
La conclusion résulte maintenant
du fait que chaque $\tau_{\leq r-i}\cp_{X_i/S}(\Lambda)$ est constructible
---~hypothèse satisfaite si $i>0$ (récurrence sur $r$) et pour $i=0$ (cas plurinodal,
traité plus bas)~--- et par la suite spectrale, comme en \ref{section-descente}, 
du fait que $\tau_{\leq r-i} \fl{\varepsilon_{i}}_{*} (K)
\giso \tau_{\leq r-i} \fl{\varepsilon_i}_{*} \tau_{\leq r-i}(K)$, ce dernier complexe étant
constructible sur $X\timesf_S S$ si $\tau_{\leq r-i}(K)$ l'est sur $X_i\timesf_S S$.

Remarquons également qu'en vertu du lemme précédent, si l'on sait
démontrer la constructibilité pour un complexe de cycles proches après
un changement de base propre et surjectif, on aura automatiquement la
constructibilité sur toute base telle que l'on ait commutation aux
changements de base.

\subsection{Le cas plurinodal et des coefficients constants ; fin de
la démonstration de \ref{constructibilite}} Soit $r$ un entier fixé tel que l'on ait constructibilité (après changement de base)
de $\tau_{\leq r-1}\cp_m\mc{F}$ pour tout morphisme $m$ comme en \ref{énoncé} et tout faisceau $\Lambda$-constructible
$\mc{F}$ sur la source de $m$. On va montrer qu'il en est de même de $\cp_f{\Lambda}$, avec $f$ plurinodal,
en procédant par récurrence sur la dimension relative. Considérons
qu'un morphisme fini est plurinodal de dimension relative nulle. 
La constructibilité dans ce cas est facile (cf. aussi \textit{infra}).

Factorisons le morphisme
plurinodal $f$ comme en \ref{dim1}. Notons $Z\sr{i}{\hra}X$ le lieu
singulier fermé de $h$, $X^0\sr{j}{\hra } X$ l'ouvert complémentaire, $\pi$
le morphisme fini $Z\ra Y$ et enfin $h^0$ le morphisme lisse $X^0\ra
Y$.  On a un triangle distingué $$ \fl{j}_! \cp_{X^0/S}(\Lambda)\ra
\cp_{X/S}(\Lambda)\ra \fl{i}_*\big(
\fl{i}^*\cp_{X/S}(\Lambda)\big)\sr{+1}{\ra},$$
et une suite exacte (où l'on note  $\cp^i$ pour $\HH^i\cp$, $i\in \ZZ$) :
$$
0\ra \fl{j}_!\cp_{X_0/S}^r\Lambda\ra \cp_{X/S}^r\Lambda\ra \fl{i}_*
\fl{i}^*\cp_{X/S}^r\Lambda\ra 0.
$$
Par hypothèse de récurrence sur la dimension relative,  nous pouvons
supposer que le $r$-tronqué de $\fl{j}_!
\cp_{X^0/S}(\Lambda)$ est constructible, car, par locale acyclicité des morphismes
lisses, $\cp_{X^0/S}\Lambda\giso
\fl{h^0}^*\cp_{Y/S}\Lambda$.
Il nous rester alors à vérifier que
$\fl{i}^*\cp_{X/S}^r(\Lambda)$ est constructible --- du moins après
changement de base. (Dorénavant, nous omettrons cette précaution
oratoire.)   
En appliquant le foncteur $h_*$, on obtient la suite exacte longue de cohomologie :

$$
0\ra \R^0 \fl{h}_* \fl{j}_! \cp_{X_0/S}^r\Lambda\ra \R^0 \fl{h}_* \cp_{X/S}^r\Lambda
\ra \R^0 \fl{h}_* \fl{i}_*
\fl{i}^* \cp_{X/S}^r\Lambda = \fl{\pi}_*\big( \fl{i}^*\cp_{X/S}^r\Lambda \big) \ra
\R^1 \fl{h}_*\fl{j}_!  \cp_{X_0/S}^r\Lambda\ra \cdots.
$$

Comme $\fl{j}_! \cp_{X_0/S}^r\Lambda$ est constructible, il en est de même des 
premiers et derniers termes (\ref{constrimagedirecte}). 
Par finitude de $\pi$, la constructibilité de $\fl{i}^*\cp_{X/S}^r\Lambda$
résultera de celle de $\fl{\pi}_*\big( \fl{i}^*\cp_{X/S}^r\Lambda \big)$.
On en tire qu'il nous suffit de montrer que le second terme, $\R^0 \fl{h}_* \cp_{X/S}^r\Lambda$
est constructible. Du lemme ci-dessous, appliqué à $K=\cp_{X/S}\Lambda$, il résulte qu'il suffit de s'assurer de la constructibilité
de  $\cp_{Y/S}^r (h_* \Lambda)=\HH^r\fl{h}_* \cp_{X/S}\Lambda$ et de $\tau_{\leq r-1}\cp_{X/S}\Lambda$ (connue par hypothèse
de récurrence sur $r$).
Or on sait que le \textit{complexe} $h_* \Lambda$ de faisceaux sur $Y$ est constructible.
Grâce aux dévissages des sections précédentes, et de l'hypothèse de double récurrence (sur l'indice $r$ et la 
dimension relative) on sait que le faisceau $\cp_{Y/S}^r (h_* \Lambda)$ est constructible.
En effet, on se ramène à la constructibilité de faisceaux $\cp^r_{Y'/S}\Lambda$
pour $Y'/S$ plurinodal de dimension relative inférieure ou égale à celle de $Y/S$ ainsi que
de faisceaux $\cp^{r-i}_m\Lambda$ pour $i>0$.

Ceci conclut la démonstration du théorème \ref{constructibilite}.

\begin{lmm2}
Soient $K\in \ob \mathsf{D}^+(X\timesf_S S,\Lambda)$ et $r\in \ZZ$. Supposons $\tau_{\leq r-1}K$ constructible
sur $X\timesf_S S$ et $\HH^r\fl{h}_*K$ constructible sur $Y\timesf_S S$.
Alors le faisceau $\R^0 \fl{h}_* \HH^r K$ est constructible sur $Y\timesf_S S$.
\end{lmm2}

\begin{proof}
On peut supposer que $\tau_{\leq r}K \ra K$ est un isomorphisme. Appliquant le foncteur
$\fl{h}_* $ au triangle distingué 
$$
\tau_{\leq r-1}K \ra K=\tau_{\leq r}K\ra \HH^r K[-r] \sr{+1}{\ra}
$$
on obtient (en passant à la cohomologie de degré $r$) la suite exacte longue :

$$
\cdots \ra \HH^r\fl{h}_*K\ra \R^0 \fl{h}_* \HH^r K \ra \HH^{r+1}\fl{h}_*\tau_{\leq r-1}K\ra \cdots
$$

Les hypothèses du lemme nous affirment que le premier terme est constructible ainsi que le troisième
compte tenu de la propreté de $h$.
\end{proof}

\section{Un exemple}\label{unexemple}

Soient $k$ un corps séparablement clos et $S$ l'hensélisé du plan $\aff_k^2$ en l'origine $o$.
Notons $e:X=\mathrm{\acute{E}cl}_o S\ra S$ l'éclatement du point fermé.
Nous allons montrer que pour $n\geq 2$ et $\Lambda=\ZZ/n$, le complexe des cycles proches $\cp_{e}\Lambda$ 
n'est pas constructible et que sa formation ne commute pas aux changements de base.
(Si l'entier $n$ est inversible sur $k$, le premier point entraîne le second compte
tenu de \ref{constructibilite}.)
Pour chaque $D$ image inverse sur $S$ d'une droite du plan passant
par l'origine, notons $\eta_D$ son point générique, $\etb_D$ un point générique géométrique et enfin 
$o_D$ le point de $X$ correspondant, localisé sur le diviseur exceptionnel.
Calculons, pour $D,D'$ deux droites comme ci-dessus,

$$(\cp_e\Lambda)_{(o_{D'},\etb_D)}=\RG(X(o_{D'})\times_{S} S(\etb_D),\Lambda).$$

La mise en défaut de la commutation aux changements de base quelconques
se voit déjà en comparant la cohomologie des \emph{tubes} de Milnor aux \emph{fibres}
de Milnor : pour $D\neq D'$, le schéma non vide $X(o_{D'})\times_{S} S(\etb_D)$
n'a pas la même cohomologie que le schéma vide $X(o_{D'})\times_{S} \etb_D$.

La restriction de $e$ à $S(\etb_D)$ étant un isomorphisme, le tube de Milnor
$ X(o_{D'})\times_{S} S(\etb_D)$ est isomorphe à $X(o_{D'})\times_{X} X(\etb_D)$,
où $\etb_D$ est ici le point géométrique de $X$ au-dessus du point $\etb_D$ de $S$.
Il résulte de \cite{Joins@Artin} (cf.~\ref{annulation}) que ce joint n'a de cohomologie
qu'en degré nul. 
Pour montrer que l'inconstructiblité, il faut comparer le $\pi_0$ de ces joints suivant les
cas $D'=D$ et $D'\neq D$. 

Si $D'\neq D$, le joint $X(o_{D'})\times_{X} X(\etb_D)$ est concentré au-dessus du point générique $\eta_X$ de $X$ :
$$X(o_{D'})\times_{X} X(\etb_D)\giso (X(o_{D'})\times_{X} X(\etb_D))\times_X \eta_X 
\isononcan X(o_{D'})_{\eta_X}\times_{\eta_X} X(\etb_D)_{\eta_X}.$$ Remarquons que les $\eta_X$-schémas
$X(o_{D'})_{\eta_X}$ sont tous isomorphes pour $D'$ variable.

Si $D=D'$, le joint correspondant est de même la somme disjointe de $X(o_{D})_{\eta_X}\times_{\eta_X} X(\etb_D)_{\eta_X}$,
qui est isomorphe à $
X(o_{D'})_{\eta_X}\times_{\eta_X} X(\etb_D)_{\eta_X}$, et d'une \emph{composante supplémentaire}
$X(o_{D})\times_{X} \etb_D \isononcan D(o_D)\times_D \etb_D=\etb_D$.

On laisse au lecteur le soin d'en déduire que les cycles évanescents 
ne sont pas constructibles au sens indiqué en \ref{constructibilite}.

Enfin on vérifie aisément que localement, la première projection
$$e':=\mathrm{\acute{E}cl}_o\aff^2_k\times_{\aff^2_k} \mathrm{\acute{E}cl}_o\aff^2_k \ra \mathrm{\acute{E}cl}_o\aff^2_k $$
est de la forme $\aff^1_k \times_k f$,
où $f$ est la projection
$$\SP(k)[X,T]/(TX)\ra \SP(k)[X].$$
En particulier, le but du morphisme $f$ est régulier de dimension $1$, ce qui entraîne la commutation
aux changements de base et la constructibilité de $\cp_f\Lambda$ et finalement les mêmes propriétés
pour $e'$. Ainsi, l'éclatement tue l'éclatement.

\section{Conjugaison des cycles évanescents à la manière de S.~Lefschetz.}

\subsection{}Soient $f:\mc{X}\ra \mc{S}$ un morphisme \textit{propre}
de schémas, $\fl{f}:\mc{X}\timesf_{\mc{S}} \mc{S} \ra
\mc{S}\timesf_{\mc{S}} \mc{S}$ le morphisme  induit, $p_1,p_2$
respectivement les première et seconde projections
$\mc{S}\timesf_{\mc{S}} \mc{S}\rra \mc{S}$  et enfin $\ell$ un nombre
premier inversible sur $S$.  On a une flèche d'adjonction $p_2^*\R
f_*\QL\ra \R\fl{f}_*\big(\cp_f\QL\big)$ associée au diagramme
commutatif

$$\xymatrix{\mc{X} \ar[d]^f \ar[r]^{\cp_f} & \mc{X}\timesf_{\mc{S}}
\mc{S} \ar[d]^{\fl{f}} \\ \mc{S} & \mc{S}\timesf_{\mc{S}} \mc{S}
\ar[l]^{p_2} } $$

C'est un isomorphisme (cf. \ref{calcfibre}) . De plus on a un
morphisme canonique de complexes de faisceaux sur
$\mc{S}\timesf_{\mc{S}} \mc{S}$, $p_1^*\R f_*\QL\ra p_2^*\R f_*\QL$
dont la fibre en un point $t\leadsto s$ s'identifie au morphisme de
spécialisation $\RG(\mc{X}_s,\QL)\ra \RG(\mc{X}_t,\QL)$.  Tout cône de
ce morphisme est isomorphe à $\R \fl{f}_*\big(\ce_f\QL\big)$.  (Le
lecteur se convaincra aisément que  le morphisme $p_1^*\R f_*\QL\ra
p_2^*\R f_*\QL$ s'obtient en appliquant le foncteur  $\R \fl{f}_*$ aux
deux premiers termes du triangle distingué sur $\mc{X}\timesf_{\mc{S}}
\mc{S}$ : $\QL\ra\cp_f(\QL)\ra\ce_f(\QL)\sr{+1}{\ra}$.)  Dans le
groupe de Grothendieck des faisceaux sur $\mc{S}\timesf_{\mc{S}}
\mc{S}$, on a donc plus suggestivement : $$[\big(\R
f_*\QL\big)_{\textrm{gén.}}]-[\big(\R f_*\QL\big)_{\textrm{spé.}}]=
[\R \fl{f}_*\big(\ce_f\QL\big)]=:[\varphi_f(\QL)].$$

\subsection{}Soit maintenant $X$ une sous-variété projective lisse
connexe de dimension $n+1$  d'un espace projectif $\PP$ sur un
corps algébriquement clos. Nous reprenons les notations et la terminologie de
\cite{weilii}, \S~4.2.3.  Soit $D\subset \PP^{\vee}$ un pinceau de
Lefschetz (non nécessairement \textit{transverse}) 
et $S=D\cap X^{\vee}$ le schéma des points $s$ de
la droite projective $D$ qui correspondent à des fibres $X\cap H_s$
singulières.  Soient $t$ un point générique géométrique de $D-S$, et
$s$ un point géométrique de $S$.

Pour chaque chemin $c$ entre $t$ et un point générique géométrique de
$D(s)$, on définit un cycle évanescent  $\pm \delta_c\in
\HH^n(X_t,\QL)([\frac{n}{2}])$. Dorénavant nous nous autoriserons
à omettre les twists à la Tate.

Comme suggéré dans \textit{loc.~cit.}, le formalisme des cycles
évanescents sur une base quelconque permet de mieux comprendre
la conjugaison des cycles évanescents (au signe près) sous l'action
du groupe $\gp(D-S,t)$, y compris dans le cas exceptionnel 
($p=2$, dimension des sections hyperplanes paire) \og sauvage \fg et/ou
des coefficients de torsion.

Notons ${X^{\vee}}_{\mathrm{bon}}$ le lieu 
de $X^{\vee}$ correspondant à une unique singularité quadratique 
ordinaire (\sga{7}{xvii}{3.2}), et supposons que $X^{\vee}$ soit une hypersurface. 
Rappelons que dans le cas non exceptionnel, 
${X^{\vee}}_{\mathrm{bon}}$ est soit vide soit le lieu lisse  ${X^{\vee}}_{\mathrm{lisse}}$
de $X^{\vee}$ (\textit{loc. cit.}~3.5). Dans le cas
exceptionnel, cet ouvert est contenu dans ${X^{\vee}}_{\mathrm{lisse}}$ mais 
l'inclusion peut être stricte.
Soit $\mc{X}$ la «variété d'incidence» des
points $(x,H)$ avec \mbox{$X\ni x\in H$} et $f$ le morphisme propre
$\mc{X}\ra \PP^{\vee}$, $(x,H)\mapsto H$. Considérons l'anneau de coefficients
$\Lambda=\ZZ/\ell^i\ZZ$, où $\ell$ est un nombre premier inversible sur $k$ et $i\in 
\NN-\{0\}$. Les résultats qui
suivent sont donc, même dans le cas non exceptionnel, légèrement plus forts que ceux
de \sga{7}{xviii}{} ($\QL$-coefficients).

D'après la théorie de Lefschetz locale, la restriction du faisceau
constructible $\varphi^n_f\Lambda$ au sous-topos 
$\fl{U}:={X^{\vee}}_{\mathrm{bon}}\timesf_{\PP^{\vee}}
\big(\PP^{\vee}-X^{\vee}\big)$ est localement
constante de rang $1$ ; les fibres ont un générateur canonique au signe près,
compatible avec les flèches de spécialisations (la démonstration, qui généralise
celle de \sga{7}{xv}{2.2}, est laissée au lecteur).
De même, $p_2^*\R f_*\Lambda=:\big(\R f_*\Lambda\big)_{\textrm{gén.}}$ est localement 
constant sur $\fl{U}$, de
fibre la cohomologie d'une section hyperplane lisse.  Le noyau $\mc{K}$ du
morphisme $\big(\R^n f_*\Lambda\big)_{\textrm{gén.}}\ra \varphi^n_f\Lambda$
est donc lisse (sur $\fl{U}$) ; sa fibre en un point $t\leadsto s$ est
isomorphe à l'orthogonal de $\pm \delta^t_s$ dans $\HH^n(X\cap
H_t,\Lambda)$. C'est aussi l'image de $p_1^*\R^n f_*\Lambda=:(\R^n
f_*\Lambda)_{\textrm{spé.}}$ dans $(\R^n f_*\Lambda)_{\textrm{gén.}}$.
L'accouplement parfait $\R^n f_*\Lambda\otimes \R^n f_*\Lambda\ra \Lambda$ sur
$\PP^{\vee}-X^{\vee}$ induit un tel accouplement sur $(\R^n
f_*\Lambda)_{\textrm{gén.}}$ par image inverse sur $\fl{U}$.  Notons
$\mathrm{Ev}$ l'orthogonal pour cet
accouplement du sous-faisceau $\mc{K}$ précédent ; c'est un 
faisceau lisse de rang $0$ ou $1$ sur
$\fl{U}$, dont la fibre en $t\leadsto s$ est canoniquement engendrée
par $\pm \delta^t_s$.  Comme le schéma ${X^{\vee}}_{\mathrm{bon}}$ est connexe, le
topos $\fl{U}$ est «connexe par arcs» : on peut relier deux points par
une chaîne finie de spécialisation/générisation de points.  Soient
maintenant deux points $s_1$ et $s_2$ de $S$ et $t$ une générisation
géométrique commune dans $D$.  Par connexité, il existe un chemin,
\cad un isomorphisme de foncteurs fibres, $\fl{g}\in 
\gp({X^{\vee}}_{\mathrm{bon}}\timesf_{\PP^{\vee}} (\PP^{\vee}-X^{\vee}) ; (s_1,t),(s_2,t))$
envoyant $\pm \delta^t_{s_1}\in \mathrm{Ev}_{(s_1,t)}$ sur $\pm
\delta^t_{s_2}\in \mathrm{Ev}_{(s_2,t)}$.  Ces éléments sont
naturellement dans $\HH^n(X_t,\Lambda)$ ; sur ce groupe, l'action de
$\gp({X^{\vee}}_{\mathrm{bon}}\timesf_{\PP^{\vee}} (\PP^{\vee}-X^{\vee}) ; (s_1,t),(s_2,t))$ 
se factorise à travers $\gp(\PP^{\vee}-X^{\vee}, t)$
par la seconde projection, donc l'image $g\in \gp(\PP^{\vee}-X^{\vee},
t)$ de $\fl{g}$ conjugue $\pm \delta^t_{s_1}$ à $\pm \delta^t_{s_2}$.
La conjugaison des cycles évanescents sous l'action de $\gp(D-S)$ 
résultera alors du fait que la représentation $\rho:\gp(\PP^{\vee}-X^{\vee},u)\ra
\Aut\big((\R^n f_*\Lambda)_u\big)$ associée à $\R^n f_*\Lambda$ a même image
que le morphisme composé $\rho\circ \big(\gp(D-S)\ra \gp(\PP^{\vee}-X^{\vee})\big)$.
Dans le cas modéré cela résulte de la surjectivité de $\gp^{\mathrm{mod}.}(D-S)
\ra \gp^{\mathrm{mod}.}(\PP^{\vee}-X^{\vee})$ (cf. \textit{loc. cit.}) ;
cela suffit donc pour conclure dans ce cas
(i.e. $p\neq 2$ ou $n$ impair)\footnote{Remarquons qu'en caractéristique positive,
le morphisme $\gp(D-D\cap H)\ra \gp(\PP-H)$ n'est \emph{jamais} surjectif si $D\neq \PP$
et $H$ est une hypersurface. Par la théorie d'Artin-Schreier, il suffit de montrer que
si $d\geq 1$, $f\in \Gamma(\PP²,\mc{O}(d))$ et $D:=V(x_0)\varsubsetneq V(f)=:H$,
il existe une fonction $g$ sur $\PP²-H$, nulle sur $D$, qui n'est pas de la forme
$h^p-h$ pour $h\in \Gamma(\PP²-H,\mc{O})$. On remarque alors que $g=\frac{x_0^d}{f}$ convient.}. \\
Nous présentons ici une démonstration de l'égalité de ces groupes de monodromie, 
due à O.~Gabber\footnote{Lettre à l'auteur, 14 mars 2005.} qui permet donc 
de traiter le cas général.
Il s'agit de montrer que si $Y\ra \PP^{\vee}-X^{\vee}$ est le revêtement galoisien
correspondant à $\mathrm{Im}(\rho)$, le schéma $Y\times_{\PP^{\vee}} D$ est connexe.
Nous allons, par un théorème de Bertini, nous ramener au cas où $D$ est la droite générique
après avoir convenablement compactifié la situation.

\begin{prp}
Soient $S=\SP(A)$ un schéma strictement local de point fermé $s$, 
$X$ un $S$-schéma plat, de présentation finie
dont la fibre spéciale est purement de dimension paire $n$ et présente 
une singularité quadratique 
ordinaire en $x$ telle que $\kappa(x)\diagup \kappa(s)$ soit radicielle.
Il existe une forme quadratique non dégénérée $\sum_{1\leq i,j\leq n} a_{i,j}X_i X_j$ 
à coefficients dans $A$ et $b,c\in A$ tels que l'hensélisé 
$X(x)$ de $X$ en $x$ soit $S$-isomorphe à l'hensélisé de 
$$S[X_0,X_1,\dots,X_n]/\big(X_0^2+b\cdot X_0+c+\sum_{1\leq i,j\leq n} a_{i,j}X_i X_j\big)$$
en le point au-dessus de $s$ dont les $n$ dernières coordonnées sont nulles.
De plus, la $S$-classe d'isomorphisme de $S[X_0]/\big(X_0^2+b\cdot X_0+c\big)$
est indépendante des choix. 
\end{prp}

\begin{proof}
Seul le dernier point est à vérifier ; le premier se trouve dans \sga{7}{xv}{1.3.2}
(dont la démonstration est valable sous réserve que $S$ soit strictement hensélien ou 
encore $\kappa(x)\diagup \kappa(s)$ triviale)\footnote{D'après O.~Gabber,
le $\SP\big(\FF_2(b)\big)$-schéma défini par l'équation $X_0^2+b+X_0(X_1^2+X_1X_2+X_2^2)$
est un contre-exemple à \textit{loc.~cit.}}.
Nous ne vérifions cette proposition que dans le cas (plus délicat) où la
caractéristique résiduelle est égale à $2$ ; nous l'appliquerons dans ce cas
uniquement.
Comme la forme quadratique $\sum_{1\leq i,j\leq n} a_{i,j}X_i X_j$ est non dégénérée,
le lieu singulier relatif est défini par l'idéal $(2X_0+b,X_1,\dots,X_n)$.
Ainsi, si l'on a une autre description de $X(x)$ décrite avec des variables 
$\star'$, on a l'égalité
$(2X_0+b,X_1,\dots,X_n)=(2X'_0+b',X'_1,\dots,X'_n)$.

Il en résulte qu'il existe des éléments $(u_{i,j})$ de $A$ tels que 
$$X_i'=\big(\sum_1^n u_{i,j} X_j\big) + v_i\cdot(2X_0+b).$$
La matrice $(u_{i,j})$ est inversible car dans la fibre spéciale
$(X_1,\dots,X_n)$ et $(X_1',\dots,X_n')$ sont des systèmes minimaux de générateurs
de l'idéal du lieu singulier relatif et $2X_0+b$ s'annule sur la fibre
spéciale ($2=0$ et $b^2-4c=0$). Finalement, 
$$(X_1',\dots,X_n')=\big(X_1-h_1\cdot(2X_0+b),\dots,X_n-h_n\cdot(2X_0+b)\big)$$
pour certains $h_i$. Ainsi, $A[X'_0]/({X'_0}^2+b'X'_0+c')\isononcan 
\mc{O}\big(X(x)\big)/(X'_1,\dots,X'_n)$ est isomorphe à 
$$\mc{O}\big(\aff^{n+1}_S(x)\big)/\big(X_i-\tilde{h}_i\cdot(2X_0+b)\,(1\leq i \leq n)\, ,
X^2_0+bX_0+c+\varphi\cdot(2X_0+b)^2\big)$$ pour certains $(\tilde{h}_i)$ et $\varphi$ dans 
$\mc{O}\big(\aff^{n+1}_S(x)\big)$. En quotientant par les $n$ premières équations,
on obtient l'anneau local hensélien $B$ de la droite affine $\aff^1_S$ de coordonnée
$X_0$ au point correspondant à $x$. On achève la démonstration en remarquant que 
si $\psi$ est une solution de $\varphi=\psi+\psi^2$, et que l'on pose
$Y_0=X_0+\psi\cdot(2X_0+b)$, on a l'égalité $Y_0^2+bY_0+c=X^2_0+bX_0+c+\varphi\cdot(2X_0+b)^2$.
\end{proof}

La proposition précédente se reformule de la façon suivante.
Soit $\Omega=\PP^{\vee}-(X^{\vee}-{X^{\vee}}_{\mathrm{bon}})$
et considérons le faisceau étale associé au préfaisceau 
{\small
$$
U\diagup \Omega\leadsto \{\text{ensemble des classes d'isomorphisme de morphismes
finis, plats de rang\ } 2, U'\ra U\ \text{étales hors de}\ X^{\vee}\}.
$$}La correspondance $\mc{X}\ra \PP^{\vee}$ définit une section de ce faisceau
sur $\Omega$. La lissité de $\mc{X}$ entraîne la lissité des schémas $U'/U$
correspondant à ces revêtements doubles. De même, si $D$ définit un pinceau
de Lefschetz, la partie correspondante de $\mc{X}$ est un éclatement lisse de la variété
$X$ si bien que $U'\times_{\PP^{\vee}} D$ est également lisse. 
Appliquons la formule de Picard-Lefschetz en le point générique de ${X^{\vee}}_{\mathrm{bon}}$.
En dimension relative paire, elle s'écrit
$$\sigma(x)-x=\pm \frac{\varepsilon(\sigma)-1}{2}\langle x,\delta\rangle \delta,$$
où le caractère $\varepsilon$ d'ordre $2$ est précisément tué par le revêtement
introduit plus haut. Ainsi, pour $U'$ comme ci-dessus et compte tenu de la pureté 
du lieu de ramification, l'image inverse du faisceau $\R^n f_*\Lambda$ sur 
$U'\times_{\PP^{\vee}} (\PP^{\vee}-X^{\vee})$ se prolonge en un faisceau lisse sur $U'$.
De plus, pour tout point géométrique $t$ de ${X^{\vee}}_{\mathrm{bon}}$, la restriction
de la représentation $\rho$ à $\gp(\PP^{\vee}(t)-X^{\vee}(t))$ est non triviale
(car $\langle \delta,\delta \rangle=\pm 2$) si bien que si $\sur{Y}$ est la normalisation
de $\PP^{\vee}$ dans $Y$, le schéma $\sur{Y}\times_{\PP^{\vee}}\PP^{\vee}(t)$ est une
union disjointe de copies des revêtements doubles distingués de $\PP^{\vee}(t)$.
(En effet, il est étale sur $\PP^{\vee}(t)-X^{\vee}(t)$, non trivial et trivialisé
par $\PP^{\vee}(t)'$.) Cela montre en particulier qu'au-dessus de $\Omega$, 
le schéma $\sur{Y}$ est lisse et, de façon semblable, que pour chaque pinceau de Lefschetz $D$,
le schéma $\sur{Y}\times_{\PP^{\vee}} D$ est lisse. Ces schémas forment
une famille propre et lisse sur l'ouvert de $\mathrm{Gr}(1,\PP^{\vee})$ correspondant
aux pinceaux de Lefschetz. La fibre générique est donc géométriquement connexe
par le théorème de Bertini (\cite{Bertini@Jouanolou}, I~6.10(3)) et il en est donc
ainsi de chaque fibre.

\subsection{Remarques finales}

Dans \sga{7}{XV}{}, P.~Deligne démontre la formule de Picard-Lefschetz
en  dimension relative impaire en utilisant un théorème de comparaison
avec  la théorie transcendante. Récemment, L.~Illusie en a donné une
démonstration algébrique (\cite{Picard@Illusie}). Il serait cependant
intéressant de disposer d'un théorème de comparaison (après
modification) --- dont le sens exact reste encore à préciser --- 
entre la cohomologie étale des fibres de Milnor (à
coefficients constants) et la cohomologie de Betti des fibres de
Milnor classiques, définies au moyen de  petites boules de rayons
«$\varepsilon,\eta$» dans le cas où $X\ra S$ est un morphisme de
variétés algébriques complexes. Le cas des singularités isolées
devrait résulter d'un argument local-global. L'auteur remercie
P.~Deligne de lui avoir fait remarquer que le cas général 
semble inconnu et L.~Illusie d'avoir attiré son attention
sur l'article \cite{RH@IKN}.

Enfin, motivé par les succès de la théorie des faisceaux pervers sur
un $S$-schéma, avec $S$ de dimension $0$ ou~$1$,  il serait naturel
d'étudier un tel formalisme dans le présent cadre ($S$ général) :
perversité des cycles proches, dualité.

\nocite{sga4} \nocite{sga41/2}


\begin{thebibliography}{GNAPGP88}

\bibitem[Art71]{Joins@Artin} {\scshape M.~Artin} -- {\og On the joins
of {H}ensel rings\fg}, \emph{Advances in Math.} \textbf{7}
(\osn{1971}), p.~282--296.

\bibitem[Art73]{sga4xv} {\scshape M.~Artin} -- {\og Morphismes
acycliques\fg}, \oldstylenums{1973}, exposé {\sc xv} dans \cite{sga4}.

\bibitem[Del72]{sga7xv}
{\scshape P.~Deligne} -- {\og La formule de {P}icard-{L}efschetz\fg},
  \oldstylenums{1972}, exposé {\sc xv} dans \cite{sga7}.

\bibitem[Del74]{Hodge3@Deligne} {\scshape P.~Deligne} -- {\og
Th\'eorie de {H}odge. {I}{I}{I}\fg}, \emph{Publications {M}ath\'ematiques 
de l'{I}.{H}.{\'E}.{S}.}
(\oldstylenums{1974}), no.~44, p.~5--77.

\bibitem[Del77]{sga41/2} {\scshape P.~Deligne \textit{et al.}} --
\emph{Cohomologie \'etale}, Springer-Verlag, \oldstylenums{1977},
Lecture Notes in Mathematics, Vol.~569.

\bibitem[Del80]{weilii} {\scshape P.~Deligne} -- {\og La conjecture de
Weil. {\scshape ii} \fg}, \emph{Publications {M}ath\'ematiques de l'{I}.{H}.{\'E}.{S}.}
(\oldstylenums{1980}), no.~52, p.~137--252.

\bibitem[dJ97]{Families@de_Jong} {\scshape A.~J. de~Jong} -- {\og
Families of curves and alterations\fg}, \emph{Ann. Inst. Fourier}
\textbf{47} (\osn{1997}), no.~2, p.~599--621.

\bibitem[GNAPGP88]{LN1335}
{\scshape F.~Guill{\'e}n, V.~Navarro~Aznar, P.~Pascual~Gainza {\normalfont
  \smfandname} F.~Puerta} -- \emph{Hyperr\'esolutions cubiques et descente
  cohomologique}, Lecture Notes in Mathematics, vol. 1335, Springer-Verlag,
  Berlin, \osn{1988}.

\bibitem[Gro67]{ega} {\scshape A.~Grothendieck} (rédigés avec la
collaboration de {\scshape J. Dieudonné}) -- {\og \'{E}l\'ements de
g\'eom\'etrie alg\'ebrique\fg}, \emph{Publications Mathématiques de
l'{I}.{H}.{\'E}.{S}.} (\oldstylenums{1960}--\oldstylenums{1967}).

\bibitem[Gro72]{sga7} {\scshape A.~Grothendieck \textit{et al.}} --
\emph{Groupes de monodromie en g\'eom\'etrie alg\'ebrique},
Springer-Verlag, \oldstylenums{1972}, S\'eminaire de g\'eom\'etrie
alg\'ebrique du Bois-Marie 1967--1969 (SGA~7~I~\&~II). Lecture Notes
in Mathematics, Vol.~288~\&~340.

\bibitem[Gro73]{sga4} {\scshape A.~Grothendieck \textit{et al.}} --
\emph{Th\'eorie des topos et cohomologie \'etale des sch\'emas},
Springer-Verlag, \oldstylenums{1972}-\oldstylenums{1973}, S\'eminaire
de g\'eom\'etrie alg\'ebrique du Bois-Marie 1963--1964 (SGA 4).
Lecture Notes in Mathematics, Vol. 269-270 \& 305.

\bibitem[Hei73]{NoetherianSpectrum@Heinzer} {\scshape W.~Heinzer} --
{\og Minimal primes of ideals and integral ring extensions\fg},
Proc. Amer. Math. Soc., vol.~40, \oldstylenums{1973}, p.~370--372.


\bibitem[Hub96]{Adic_spaces@Huber} {\scshape R.~Huber} --
\emph{\'{E}tale cohomology of rigid analytic varieties and adic
spaces}, Aspects of Mathematics, Friedr. Vieweg \& Sohn, Braunschweig,
1996.

\bibitem[Ill02]{Picard@Illusie} {\scshape L.~Illusie} -- {\og Sur la
formule de {P}icard-{L}efschetz\fg}, Algebraic geometry 2000, Azumino
(Hotaka), Adv. Stud. Pure Math., vol.~36, Math. Soc. Japan, T\=oky\=o,
\osn{2002}, p.~249--268.

\bibitem[IKN05]{RH@IKN}
{\scshape L.~Illusie, Kat\=o~K. {\normalfont \smfandname} Nakayama~K.} -- {\og
  Quasi-unipotent logarithmic {R}iemann-{H}ilbert correspondences\fg}, \emph{J.
  Math. Sci. Univ. T\=oky\=o} \textbf{12} (\osn{2005}), no.~1, p.~1--66.

\bibitem[Jou83]{Bertini@Jouanolou} {\scshape J.-P.~Jouanolou} --
\emph{{T}h\'eor\`emes de {B}ertini et applications}, 
Progress in Mathematics, Birkh\"auser, Boston, \osn{1983}.

\bibitem[Kat72]{sga7xvii}
{\scshape N.~M. Katz} -- {\og Pinceaux de {L}efschetz : th\'eor\`emes
  d'existence\fg}, \oldstylenums{1972}, exposé {\sc xvii} dans \cite{sga7}.

\bibitem[Lau81]{Laumon-semi_cont} {\scshape G.~Laumon} -- {\og
Semi-continuit\'e du conducteur de {S}wan (d'apr\`es
{P}. {D}eligne)\fg}, Ast\'erisque, vol.~83, Soc. Math. France, Paris,
\oldstylenums{1981}, p.~173--219.

\bibitem[Lau83]{Vanishing@Laumon} {\scshape G.~Laumon} -- {\og
Vanishing cycles over a base of dimension {$\geq 1$}\fg}, Algebraic
geometry (T\=oky\=o/Ky\=oto, 1982), Lecture Notes in Math., vol. 1016,
Springer, Berlin, \osn{1983}, p.~143--150.

\bibitem[L{\"u}t93]{Nagata@Lutkebohmert}
{\scshape W.~L{\"u}tkebohmert} -- {\og On compactification of schemes\fg},
  \emph{Manuscripta Math.} \textbf{80} (\osn{1993}), no.~1, p.~95--111.

\bibitem[MV00]{Proper@Moerdijk} {\scshape I.~Moerdijk {\normalfont
\smfandname} J.~J.~C. Vermeulen} -- {\og Proper maps of toposes\fg},
\emph{Mem. Amer. Math. Soc.} \textbf{148} (\osn{2000}), no.~705.

\bibitem[OP68]{noetherian@Ohm} {\scshape J.~Ohm {\normalfont
\smfandname} R.~L. Pendleton} -- {\og Rings with noetherian
spectrum\fg}, \emph{Duke Math. J.} \textbf{35} (\osn{1968}),
p.~631--639.

\bibitem[RG71]{Platification@Raynaud} {\scshape M.~Raynaud
{\normalfont \smfandname} L.~Gruson} -- {\og Crit\`eres de platitude
et de projectivit\'e. {T}echniques de «platification» d'un module\fg},
\emph{Invent. Math.} \textbf{13} (\osn{1971}), p.~1--89.

\bibitem[Sab83]{cycles_evanescents@Sabbah} {\scshape C.~Sabbah} --
{\og Morphismes analytiques stratifi\'es sans \'eclatement et cycles
\'evanescents\fg}, Analyse et topologie sur les espaces singuliers,
(Luminy, 1981), Ast\'erisque, vol. 101, Soc. Math.  France, Paris,
\oldstylenums{1983}, p.~286--319.

\end{thebibliography}
\end{document}